\documentstyle[amssymb,amsfonts]{amsart}

\newenvironment{proof}{\noindent {\bf Proof} }{\endprf\par}
\def \endprf{\hfill  {\vrule height6pt width6pt depth0pt}\medskip}
\def\emph#1{{\it #1}}
\def\textbf#1{{\bf #1}}

\newcommand{\bea}{\begin{eqnarray}}
\newcommand{\eea}{\end{eqnarray}}
\def\beaa{\begin{eqnarray*}}
\def\eeaa{\end{eqnarray*}}

\def\ba{\begin{array}}
\def\ea{\end{array}}
\def\be#1{\begin{equation} \label{#1}}

\newcommand{\nn}{\nonumber}
\def\medn{\medskip\noindent}
\def\rrrr{{\Bbb R}}
\def\rr{{\bf R}}

\def\nn{\nonumber}
\def\stu{{S_{t,u}}}

\def\Sbb{\underline{S}}

\newcommand{\lapp}{\mbox{$\bigtriangleup  \mkern-13mu / \,$}}

\def\Lie{{\cal L}}
\def\tr{\mbox{tr}}

\def\gg{{\bf g}}
\def\mm{{\bf m}}
\def\ns{{\bf n}}
\def\vs{{\bf v}}
\def\Ext{\mbox{Ext\,}_t}
\def\Hb{{\bold H}}

\newcommand{\half}{\frac 12}

\def\Qb{\bar{Q}}
\def\f12{\frac 12}
\def\half{\frac 12}

\def\bk{\bar{k}}

\def\pit{\bar{\pi}}
\newcommand{\divv}{\mbox{div}\mkern-19mu /\,\,\,\,}

\def\a{\alpha}
\def\b{\beta}
\def\ga{\gamma}
\def\Ga{\Gamma}
\def\de{\delta}

\def\ep{\epsilon}
\def\eps{\epsilon}
\def\la{\lambda}
\def\La{\Lambda}
\def\si{\sigma}
\def\Si{\Sigma}
\def\om{\omega}
\def\Om{\Omega}

\def\th{\theta}
\def\ze{\zeta}
\def\nab{\nabla}
\def\varep{\varepsilon}

\newcommand{\trchb}{\tr \chib}
\newcommand{\chih}{\hat{\chi}}
\newcommand{\chib}{\underline{\chi}}
\newcommand{\xib}{\underline{\xi}}
\newcommand{\etab}{\underline{\eta}}

\def\f14{\frac{1}{4}}
\def\dd{{\bf D}}

\newcommand{\les}{\lesssim}
\newcommand{\ges}{\gtrsim}

\newcommand{\pik}{\,\,^{(K)}\pi}
\newcommand{\pikt}{\,\,^{(K)}\bar{\pi}}

\renewcommand{\pit}{\bar{\pi}}

\newcommand{\Tr}{\text{Tr}\,}
\newcommand{\QQ}{\mathcal{Q}}
\newcommand{\HH}{{\mathcal H}}

\newcommand{\DD}{{\mathcal D}}
\newcommand{\EE}{{\mathcal E}}
\newcommand{\EEi}{\EE^{(i)}}
\newcommand{\EEe}{\EE^{(e)}}

\newcommand{\TT}{{\cal T}}

\renewcommand\div{\mbox{div}}

\renewcommand\Sbb{\underline{S}}

\def\ub{\underline{u}}

\def\il{\int}

\def\Lb{\underline{L}}

\def\QQ{{\cal Q}}
\def\ric{\mbox{ \bf Ric}}

\def\pr{\partial}
\def\f12{\frac 1 2}
\renewcommand{\chib}{\underline{\chi}}
\def\chih{\hat{\chi}}
\def\trch{\mbox{tr}\chi}

\newcommand{\nabb}{\mbox{$\nabla \mkern-13mu /$\,}}

\begin{document}
\theoremstyle{plain}
  \newtheorem{theorem}[subsection]{Theorem}
  \newtheorem{conjecture}[subsection]{Conjecture}
  \newtheorem{proposition}[subsection]{Proposition}
  \newtheorem{lemma}[subsection]{Lemma}
  \newtheorem{corollary}[subsection]{Corollary}

\theoremstyle{remark}
  \newtheorem{remark}[subsection]{Remark}
  \newtheorem{remarks}[subsection]{Remarks}

\theoremstyle{definition}
  \newtheorem{definition}[subsection]{Definition}

\include{psfig}

\title[Nonlinear Wave Equations ]{Rough solutions of  the Einstein-Vacuum equations }
\author{Sergiu Klainerman}
\address{Department of Mathematics, Princeton University,
 Princeton NJ 08544}
\email{ seri@@math.princeton.edu}

\author{Igor Rodnianski}
\address{Department of Mathematics, Princeton University, 
Princeton NJ 08544}
\email{ irod@@math.princeton.edu}

\vspace{-0.3in}
\begin{abstract}
This is the first  in a series of papers in which 
we initiate the study of very rough  solutions  to the initial value problem
for the Einstein vacuum equations expressed relative
 to  wave coordinates. By very rough we mean solutions
which cannot be constructed by the classical techniques of energy estimates 
and Sobolev inequalities. Following \cite{Kl-Rod} we develop new 
analytic methods based on Strichartz type inequalities which results in a gain  of
half a derivative relative to the classical result.  Our methods blend paradifferential
techniques with a geometric approach to the derivation of decay estimates. The
 latter  allows
us to take full  advantage of the specific structure of the Einstein equations.

\end{abstract}
\maketitle

\section{Introduction}
We consider the Einstein Vacuum equations,
\be{I1}
\rr_{\a\b}(\gg)=0
\end{equation}
where $\gg$ is a four dimensional Lorentz metric and $\rr_{\a\b}$
its Ricci curvature tensor. In wave coordinates $x^\a$,
\be{I2}
\square_{\gg} x^\a =\frac{1}{|\gg|}\pr_\mu(\gg^{\mu\nu}|\gg|\pr_\nu)x^\a= 0,
\end{equation}
the Einstein vacuum equations take the reduced form, see \cite{Bruhat},
\cite{HKM}.
\be{I3}
\gg^{\a\b}\pr_\a\pr_\b  \gg_{\mu\nu}=N_{\mu\nu}(\gg,\pr \gg)
\end{equation}
with $N$ quadratic in the first derivatives $\pr\gg$ of the metric.
We consider the initial value problem along
the spacelike  hyperplane  $\Sigma$  given by $t=x^0=0$,
\be{I4}
\nabla \gg_{\a\b}(0)\in  H^{s-1}(\Sigma)\,,\quad  \pr_t \gg_{\a\b}(0)\in 
H^{s-1}(\Sigma)
\end{equation}
with $\nabla$ denoting the gradient with respect to the
 space coordinates $x^i$, $i=1,2,3$ and  $H^s$ the standard 
Sobolev spaces. We also assume that $\gg_{\a\b}(0)$ is
a continuous Lorentz metric and 
\be{limit}
\sup_{|x|=r}|\gg_{\a\b}(0)-\mm_{\a\b}|\longrightarrow 0\quad  \mbox{as}  \quad
r\longrightarrow \infty,
\end{equation}
where $|x|=(\sum_{i=1}^3 |x^i|^2)^{\frac{1}{2}}$ and $\mm_{\a\b}$
the Minkowski metric.

The following local existence and uniqueness
 result (well posedness) is well known (see \cite{HKM}
and the previous result of Ch. Bruhat \cite{Bruhat} for $s\ge 4$.) 

\begin{theorem} Considered the reduced equation \eqref{I3}
subject to the initial conditions  \eqref{I4} and \eqref{limit} for some $s>5/2$.
Then there exists a time interval $[0,T]$ and unique( Lorentz metric)
solution
 $\gg\in C^0([0,T]\times \rrrr^3)$,  $\partial \gg_{\mu\nu}\in C^0([0,T];
H^{s-1})$ with $T$ depending only on  the size  of the norm $\|\pr
\gg_{\mu\nu}(0)\|_{H^{s-1}}$. In addition  condition \eqref{limit}
remains true on any spacelike hypersurface $\Si_t$, i.e. any
level hypersurface of the time function $t=x^0$.  

\label{It1}
\end{theorem}

We  establish  a significant improvement of this  result bearing on the  
issue of  minimal regularity of the initial conditions:

\medn
{\bf Main Theorem}\quad  Consider a  classical solution 
 of the  equations \eqref{I3}
for which  \eqref{I1} also holds\footnote{In other words for
any solution of the reduced equations \eqref{I3} whose initial data satisfy
 the constraint equations, see \cite{Bruhat} or \cite{HKM}. The fact that our
solutions verify \eqref{I1} plays a
fundamental role in our analysis.}. We show\footnote{We assume however
 that $T$ stays sufficiently small, e.g. $T\le 1$.  This a purely technical
 assumption which one should be able to remove.} that  the
time
$T$ of existence depends in fact only on the size of the norm $\|\pr
\gg_{\mu\nu}(0)\|_{H^{s-1}}$,
 for any fixed   $s>2$. 

\begin{remark} Theorem \ref{It1} implies the  classical local  existence result of
\cite{HKM}  for asymptotically flat initial data sets $\Sigma, g, k$ with 
$\nab g, k\in H^{s-1}(\Si)$ and $s>\frac{5}{2}$, relative  to a fixed 
system of coordinates.  Uniqueness can be
proved for additional regularity  $s>1+\frac{5}{2}$. 
 We recall that an initial data set $(\Si, g, k)$ 
 consists  of
a three dimensional complete Riemannian manifold 
$(\Si, g)$, a 2-covariant symmetric tensor $k$ on
$\Si$ verifying the constraint equations:
$$
\nabla^j k_{ij}-\nabla_i \,\tr k=0
$$
$$R-|k|^2+(\tr k)^2=0$$
where $\nabla$ is the covariant derivative, $R$
the scalar curvature of $(\Si, g)$.
  An initial data set
is said to be asymptotically flat (AF) if there exists a system of coordinates
$(x^1,x^2,x^3)$
 defined in a neighborhood of infinity\footnote{We assume, for
simplicity, that $\Si$  has only one end. A neighborhood
of infinity means the complement of a sufficiently large compact set on  $\Si$.} on $\Si$
relative to which the metric
$g$ approaches the Euclidean metric  and $k$  approaches  zero\footnote{ Because of
the constraint equations 
the 
asymptotic  behavior cannot be arbitrarily prescribed. A precise definition of
asymptotic flatness   has to involve the ADM  mass of $(\Si, g)$. Taking the mass into
account we write   $g_{ij}=(1+\frac{2M}{r})\de_{ij} +o(r^{-1})$  as 
$r=\sqrt{(x^1)^2+(x^2)^2+(x_3)^2}$ $\rightarrow \infty$. . According to the positive
mass theorem $M\ge 0$ and  $M=0$ implies that the initial data set is flat. Because
of the mass term  we cannot assume that $g-e\in L^2(\Si)$, with $e$  the $3$D Euclidean
 metric.}
\end{remark}

\begin{remark} The Main Theorem  ought to imply existence and uniqueness\footnote{
Properly speaking uniqueness holds, with $s>2$, only for the reduced equations.
Uniqueness for the actual Einstein equations requires one more derivative,
see \cite{HKM}.} for initial conditions with $H^{s}$, $s>2$, regularity. To
achieve this we only need to approximate a given $H^s$ initial data set( i.e.
$\nabla
g\in H^{s-1}(\Si), k\in H^{s-1}(\Si)$, $s>2$ ) for 
the Einstein vacuum equations by classical  initial data sets, i.e. $H^{s'}$
 data sets  with 
$s'>\frac{5}{2}$, for which theorem \ref{It1} holds. The Main Theorem allows
us to pass to the limit and derive existence of solutions for the given, rough,
initial data set. We don't know however if such an approximation result 
for the constraint equations exists in the literature.
\end{remark}

For convenience we shall also write the reduced equations
\eqref{I3} in the form
\be{I3'}
\gg^{\a\b}\pr_\a\pr_\b\phi=N(\phi,\pr \phi)
\end{equation}
where $\phi=(\gg_{\mu\nu})$, $N=N_{\mu\nu}$
and $\gg^{\a\b}=\gg^{\a\b}(\phi)$.

Expressed relative to the wave coordinates $x^\a$
the spacetime metric $\gg$   takes the form:

\be{I5}
 \gg = - \ns^2 dt^2 + g_{ij} (dx^i + \vs^i dt) (dx^j + \vs^j dt)
\end{equation}
where $g_{ij}$ is a Riemannian  metric on the
slices $\Sigma_t$, given by the level
hypersurfaces of the time function $t=x^0$,    $\ns$ is the  lapse
function of the time  foliation, and
$\vs$ is a vector-valued shift function.
The components of the inverse  metric  $ \gg^{\a\b} $
 can be
found as follows:
$$
 \gg^{00} = - \ns^{-2}, \quad  \gg^{0i} = \ns^{-2} \vs^i, \quad
 \gg^{ij} = g^{ij} - \ns^{-2} \vs^i \vs^j.
$$
In view of the Lorentzian character of $\gg$ and the spacelike character
of the hypersurfaces $\Si_t$,
\be{timelike}
c|\xi|^2\le \gg_{ij} \xi^i\xi^j\le c^{-1}|\xi|^2, \qquad c \le \ns^2-
|\vs|_{g}^2
\end{equation}
for some $c>0$.

The classical local existence result for systems of wave equations
of type \eqref{I3'} is based on energy estimates and the standard
$ H^s\subset L^\infty$ Sobolev inequality. Indeed using energy estimates and simple 
commutation inequalities one can show that,
\be{Ienergy}\|\pr \phi(t)\|_{H^{s-1}}\le E\|\pr \phi(0)\|_{H^{s-1}}
\end{equation}
with a constant $E$,
\be{Iintegr}
E=\exp\bigg(C\int_0^t\|\pr\phi(\tau)\|_{L_x^\infty}d\tau\bigg)
\end{equation}
By the classical Sobolev inequality,
$$E\le \exp\bigg(Ct\sup_{0\le \tau\le t}\|\pr\phi(\tau)\|_{H^{s-1}}d\tau\bigg)
$$
provided that $s> \frac{5}{2}$. The classical local existence result follows
by combining this last estimate, for  a small time
interval,   with  the  energy estimates
\eqref{Ienergy}.

 This scheme is very wasteful. To do better one would like
to take advantage of the mixed $L_t^1L_x^\infty$ norm appearing on the right hand
side of \eqref{Iintegr}. Unfortunately there are no good estimates for such norms
even when $\phi$ is simply a solution of the standard wave equation
\be{standardwave}
\square \phi=0
\end{equation}
in Minkowski space.
There exist however improved regularity estimates
for solutions of \eqref{standardwave} in the mixed $L_t^2L_x^\infty$ norm .
More precisely, if $\phi$ is a solution of \eqref{standardwave} and $\eps>0$
 arbitrarily small,
\be{Istrich}
\|\pr\phi\|_{L_t^2L_x^\infty([0,T]\times \rrrr^3)}\le
CT^{\eps}\|\pr\phi(0)\|_{H^{1+\eps}}.
\end{equation}
Based on this fact it was reasonable to hope that
one can improve the Sobolev exponent in the classical local existence theorem from
$s>\frac{5}{2}$  to $s>2$. This can be easily done  for solutions of
semilinear equations, see \cite{Po-Si}. In the quasilinear case, however, the situation
is far more difficult. One can no longer rely on  the Strichartz 
inequality \eqref{Istrich} for the flat D'Alembertian in 
\eqref{standardwave}; we need   instead  its extension  to the  operator
$\gg^{\a\b}\pr_\a\pr_\b$  appearing in \eqref{I3'}. Moreover, since
the metric $\gg^{\a\b}$ depends on the solution  $\phi$, it can have
only as much regularity as $\phi$ itself. This means that we have
to confront the issue of proving Strichartz estimates for wave
operators $g^{\a\b}\pr_\a\pr_\b$ with \textit{ very rough coefficients} $\gg^{\a\b}$.
This issue was  recently addressed in the pioneering works of  Smith\cite{Sm},
 Bahouri-Chemin \cite{Ba-Ch1}, \cite{Ba-Ch2}  and Tataru \cite{Ta1}, \cite{Ta2}, we refer to
the introduction in
\cite{Kl} and
\cite{Kl-Rod} for a more thorough discussion of their important  contributions.

The results of Bahouri-Chemin and Tataru are based on establishing 
a Strichartz type inequality, \textit{with a loss}, for wave operators 
with  \textit{very rough 
coefficients}\footnote{The  derivatives of the  coefficients
$\gg$  are required to be  bounded in  $L_t^\infty H_x^{s-1}$ 
and $L_t^2 L_x^\infty$  norms, with $s$ compatible  with the regularity 
 required  on the 
right hand side of the Strichartz inequality one wants  to  prove.  }.
The optimal result\footnote{ Recently Smith-Tataru \cite{Sm-Ta}
have shown that the result of Tataru is indeed sharp.} in this regard,
  due to Tataru, see \cite{Ta2}, requires a loss of  $\si=\frac{1}{6}$. 
This leads to a proof of local well posedness for systems of type 
\eqref{I3'} with $s>2+\frac{1}{6}$. 

To do better than that one needs
to take into account the nonlinear structure of the equations. In 
\cite{Kl-Rod} we were able to improve the result of Tataru by 
taking into account not only the  expected regularity properties of 
the coefficients $g^{\a\b}$ in \eqref{I3'} but also the fact that they
are themselves solutions to a similar system of equations.  This
allowed us  to improve the  exponent $s$,
needed in the proof of well posedness of equations of type\footnote{The result in
\cite{Kl-Rod} applies to general equations of type \eqref{I3'} not necessarily tied
to \eqref{I1}. In \cite{Kl-Rod} we have also made the simplifying assumptions $\ns=1$
 and  $\vs=0$.}
\eqref{I3'}, to
$s> 2+\frac{2-\sqrt{3}}{2}$. Our approach was based on a combination of the paradifferential
calculus ideas, initiated  in  \cite{Ba-Ch1} and \cite{Ta2}, with a geometric treatment of
the  actual equations  introduced  in \cite{Kl}. The main improvement was due to a gain of
 conormal differentiability for solutions to the Eikonal equations
\be{firsteikonal}
H^{\a\b}\pr_\a u\pr_\b u=0
\end{equation}  where the background metric  $H$
is a properly  microlocalized  and rescaled version 
of the metric $\gg^{\a\b}$ in \eqref{I3'}.
That gain could be traced down to the fact
that a certain component of the Ricci curvature of $H$
has a special form. More precisely denoting by $L'$  the null
geodesic vectorfield associated to $u$,  $L'=-H^{\a\b}\pr_\b u \pr_\a$, and 
rescaling it
 in an appropriate fashion\footnote{such $<L, T>_{H}=1$ with $T$ is the unit normal
to the level hypersurfaces  $\Si_t$ associated to   the time function $t$,},  $L=bL'$,  
we found that the    $\rr_{LL}=${\bf Ric}$(H)(L,L)$, verifies the remarkable
identity:
\be{remident}
\rr_{LL}=L(z)-\frac 12 L^\mu L^\nu (H^{\a\b}\pr_\a\pr_\b H_{\mu\nu})+{\mbox e}
\end{equation}
where $z\le O(|\pr H|)$ and $e\le O (|\pr H|^2)$. Thus,  apart from $L(z)$ which
is to be integrated along the null geodesic flow generated by $L$,  the only
terms which depend of the second derivatives of $H$ appear in $H^{\a\b}\pr_\a\pr_\b H$ and can
therefore be eliminated with the help of  the equations \eqref{I3'}.

 In this paper we develop the ideas of \cite{Kl-Rod} further by taking full advantage of the 
Einstein equations \eqref{I1}
 in wave coordinates \eqref{I3'}. An important aspect of our analysis here
is that the  term $L(z)$ appearing on the right hand side of \eqref{remident}
vanishes identically. We make use of both the vanishing of the Ricci curvature
of $\gg$ and the wave coordinate condition \eqref{I2}. The other important
new features are  the use of energy estimates along the null hypersurfaces 
generated by the optical function $u$ and a more efficient use of the conormal
properties of the null structure equations. 

 Our work  is divided in three parts. In this paper we give all the details
in the proof of the  Main Theorem with the exception of those results which
concern the asymptotic properties of the Ricci coefficients( the Asymptotics Theorem),
the   isoperimetric  and trace inequalities on 2-surfaces.
We   give  precise  statements of these results  in section 4.
Our second paper \cite{Einst2} is dedicated to the proof of the
Asymptotics Theorem. The isoperimetric and trace inequalities
together with some other results needed in \cite{Einst2} are 
proved in our third paper \cite{Einst3}.

We strongly believe that the result of our main theorem is not sharp.
The critical  Sobolev exponent  for the Einstein equations is $s_c=\frac{3}{2}$.
 A proof of well posedness for $s=s_c$ will provide a much stronger version
of the \textit{global stability of Minkowski space} than that of \cite{Ch-Kl}.
This is completely out of reach at the present time. A more reasonable goal,
 at the present time, is to prove the \textit{ $L^2$- curvature conjecture}, see 
 \cite{Kl1}, corresponding to the exponent $s=2$. 

\section{reduction to decay estimates}
The proof of the main theorem can be reduced to a microlocal 
decay estimate. The reduction is standard\footnote{see
\cite{Kl-Rod} and the references therein}; we quickly review
here the main steps. The precise statements and their proofs are given
in section 8.

\begin{itemize}
\item \textit{ Energy estimates}
 
Assuming that $\phi$
is a  solution\footnote{i.e. a classical
solution according to theorem  \ref{It1}.} of \eqref{I3'} on  $[0,T]\times
\rrrr^3$
 we have the
apriori energy estimate:

\begin{equation}
\|\pr \phi\|_{L^\infty_{[0,T]}\dot H^{s-1}} \le 
C\|\pr \phi(0)\|_{\dot H^{s-1}}
\label{ens}
\end{equation}
with a constant $C$ depending only on
$\|\phi\|_{L^\infty_{[0,T]}L_x^{\infty}}$
 and $\|\pr\phi\|_{L^1_{[0,T]} L^\infty_x}$.

\item\textit{ Strichartz estimate} 
To prove  our Main Theorem  we need, in addition to \eqref{ens}
 an estimate of the form:
$$\|\pr \phi\|_{L^1_{[0,T]}L_x^\infty}\le C\|\pr\phi(0)\|_{H^{s-1}}$$
for any $s>2$. We accomplish it  by establishing a Strichartz type
 inequality of the form,
\be{Istrich100} 
\|\pr \phi\|_{L^2_{[0,T]}L_x^\infty}\le C\|\pr\phi(0)\|_{H^{1+\ga}}
\end{equation}
with any fixed $\ga>0$. We achieve this with the help of a bootstrap
argument. More precisely we make the assumption

\be{bootstrap} 
\|\pr\phi\|_{L^\infty_{[0,T]} H^{1+\ga}} +\|\pr \phi\|_{L^2_{[0,T]}
L^\infty_x}\le B_0,
\end{equation}
and use it to prove the 
  better estimate;

\be{mainestimate}
\|\pr \phi\|_{L^2_{[0,T]}L^\infty_x}\le C(B_0)\, 
T^{\de}.
\end{equation}
for some $\de>0$. Thus, for sufficiently small $T>0$, we find that
\eqref{Istrich100} holds true.
\item \textit{Proof of the Main Theorem}

This can be done easily by combining the energy estimates with 
the Strichartz estimate stated above.

\item \textit{Dyadic  Strichartz  Estimate}

 The proof of the Strichartz estimate
can be reduced to a dyadic version for each
$\phi^\la=P_{\la}\phi$, $\la$ sufficiently large\footnote{The low frequencies
 are much easier to treat.},  where
$P_\la$ is the 
 Littlewood-Paley projection on the space frequencies of size $\la\in 2^{\Bbb Z}$.
$$
\|\pr \phi^\la\|_{L^2_{[0,T]}L^\infty_x}\le C(B_0)\,c_\la 
T^{\de}  
\|\pr \phi\|_{H^{1+\ga}},
$$
with  $\sum_\la c_\la\le 1$.
\item
\textit{Dyadic linearization and time restriction} 

Consider the new metric  $\,\gg_{<\la}= P_{<\la}\gg=
\sum_{\mu\le 2^{-M_0} \la} P_{\mu} \gg\,$, for some sufficiently large constant 
$M_0>0$, restricted to a
subinterval
$I$ of $[0,T]$ of size $|I|\approx T\la^{-8\ep_0}$ with $\ep_0>0$ fixed
such that $\ga>5\ep_0$. Without loss of generality\footnote{In view of the translation
invariance of our estimates.} we can assume that $I=[0,\bar{T}]$.
Using an appropriate( now standard, see \cite{Ba-Ch1}, \cite{Ta2}, \cite{Kl}, 
\cite{Kl-Rod})  paradifferential linearization
together with the Duhamel principle we can reduce the proof of the dyadic
Strichartz estimate mentioned above to a homogeneous Strichartz estimate
for  the equation 
$$\gg_{<\la}^{\a\b}\pr_\a\pr_\b \psi=0,$$
with initial conditions at $t=0$ verifying,
$$(2^{-10}\la)^m\le \|\nab^m \pr \psi(0)\|_{L^2_x}\le (2^{10}\la)^m \|\pr \psi(0)\|_{L^2_x}.$$
There exists a sufficiently small $\de> 0$,  $5\ep_0 +\de<\ga$, such that
\be{PlaStr}
\|P_\la\,\pr \psi\|_{L^2_{I} L^\infty_x}\le C(B_0)\,  
\bar{T}^\de  \|\pr \psi(0)\|_{\dot H^{1+\de}}
\end{equation}
\item \textit{Rescaling} 

Introduce the rescaled metric\footnote{$H_{(\la)}$ is a Lorentz metric 
for $\la\ge \La$ with $\La$ sufficiently large. See the discussion 
following \eqref{h} in section 8.}
$$ H_{(\la)}(t,x)=\gg_{<\la}(\la^{-1}t, \la^{-1}x)$$
and consider the rescaled  equation
$$H_{(\la)}^{\a\b}\pr_\a\pr_\b\psi=0$$
 in  the region $[0,t_{*}]\times \rrrr^3$ with $t_*\le \la^{1-8\ep_0}$
Then, with $P=P_1$, 
$$\|P\,\pr \psi\|_{L^2_{I} L^\infty_x}\le C(B_0)\,  
t_*^\de  \|\pr \psi(0)\|_{L^2}
$$
would imply the estimate \eqref{PlaStr}.
\item
\textit{Reduction to an $L^1-L^\infty$  decay estimate}

The standard way to prove a Strichartz inequality of the type
discussed above is to reduce it, by a $TT^*$ type argument,
to an $L^1-L^\infty$  dispersive type inequality.  The inequality we need,
concerning the initial value problem  
$$
 \square_{H_{(\la)}}\psi= \frac{1}{\sqrt{|H_{(\la)}|}}\pr_\a\big(H^{\a\b}_{(\la)}
\sqrt{|H_{(\la)}|}\,\pr_\b \psi\big)=0, 
$$
with data at $t=t_0$  
has the form,
$$\|P\, \pr \psi (t)\|_{L^\infty_x}\le \,  
C(B_0)\bigg (\frac 1{(1+|t-t_0|)^{1-\de}} + d(t)\bigg)
\sum_{k=0}^m \|\nab^k \pr \psi(t_0)\|_{L^1_x}$$
for some integer $m\ge 0$.

\item
\textit{ Final reduction to a localized   $L^2-L^\infty$  decay estimate}

We state this as the following theorem:
\end{itemize}
\begin{theorem}
Let $\psi$ be a solution of the equation,
\be{Idecay}\square_{H_{(\la)}}\psi=0
\end{equation}
on the time interval $[0,t_*]$ with $t_*\le \la^{1-8\ep_0}$
 Assume that the  initial data is given at $t=t_0\in [0,t_*]$,
supported in the ball  $B_{\frac 1 2}(0)$ of radius $\f12$ centered
at the origin. We fix a  large constant $\La>0$ and consider only the frequencies
$\la\ge \La$.  There exists a function $d(t)$, with 
$t_*^{\frac 1q} \|d\|_{L^q([0, t_*])}\le 1$ for some $q>2$ sufficiently close to $2$, an
arbitrarily  small $\de>0$ and a sufficiently large integer $m>0$ such
that for all $t\in[0,t_*]$, 
\begin{equation}
\|P\, \pr \psi (t)\|_{L^\infty_x}\le \,  
C(B_0)\bigg (\frac 1{(1+|t-t_0|)^{1-\de}} + d(t)\bigg)
\sum_{k=0}^m \|\nab^k \pr\psi(t_0)\|_{L^2_x}.
\label{Idecayestimate}
\end{equation}
\label{L2decay}
 \end{theorem}
\begin{remark}
\label{lesremark}
In view of the proof of the Main Theorem
 presented above, which relies on the 
final  estimate \eqref{mainestimate},  we can in
what follows treat the bootstrap constant $B_0$ as a universal constant 
and  bury
the dependence on it
  in the notation $\les$ we introduce below.   
\end{remark}
\begin{definition}
\label{lesdefinition} We use the notation $A\les B$ to express
the inequality $A\le C B$ with a universal constant, which
may depend on $B_0$ and various other  parameters
 depending only   on $B_0$  introduced
in the proof.

\end{definition}

The proof of theorem \ref{L2decay} relies on a
generalized  Morawetz type energy estimate   which will be presented in the next section.
 We shall in fact  construct  a vectorfield, analogous to the Morawetz vectorfield
in  the Minkowski space,  which depends heavily on the ``background metric'' 
$H=H_{(\la)}$. In the next  proposition we   display most  of   the main properties
 of the metric
$H$ which will be  used in the following section.
\begin{proposition}[Background estimates] 
\label{propH}
 We fix the region
$[0,t_*]\times \rrrr^3 $, with $t_*\le \la^{1-8\ep_0}$, where the original
 Einstein metric\footnote{recall that in fact  $\gg$ is  $\phi^{-1}$. Thus, in view
of the non degenerate Lorentzian character of $\gg$ the bootstrap 
assumption for $\phi$ reads as an assumption for $\gg$.} $\gg=\gg(\phi)$ verifies
the bootstrap assumption \eqref{bootstrap}.  The metric
 \be{Hla}
H(t,x)=H_{(\la)}(t,x)=P_{<\la}\gg(\la^{-1}t, \la^{-1} x)
\end{equation}
can be decomposed 
relative to our spacetime  coordinates.
\be{decompHH}
H=-n^2dt^2+h_{ij}(dx^i+v^i dt)\otimes(dx^j+v^j dt)
\end{equation}
where $n$ and $v$ are related to ${\bf n}$, ${\bf v}$
according to the rule \eqref{Hla}. The metric components $n, v$, and $h$
satisfy the conditions
\be{psithH}
c|\xi|^2\le h_{ij}\xi^i\xi^j \le c^{-1}|\xi|^2,\quad  n^2-|v|^2_h\ge c>0, \quad
|n|,|v|\le c^{-1}
\end{equation}
In addition, the derivatives of the metric $H$ verify the following:
\begin{align}
&\|\pr^{1+m}H\|_{L^1_{[0,t_*]} L_x^\infty}\les \la^{-8\eps_0},
\label{asH1}\\
&\|\pr^{1+m}H\|_{L^2_{[0,t_*]} L_x^\infty}\les
\la^{-\frac 12-4\eps_0},
\label{asH2}\\
&\|\pr^{1+m} H\|_{L^\infty_{[0,t_*]} L_x^\infty}\les 
\la^{-\f12 - 4\ep_0},
\label{asH3}\\
&\|\nab^{\frac 12+m}(\pr H)\|_{L^\infty_{[0,t_*]} L_x^2}\les \la^{-m}\qquad
\mbox{for}\quad  -\f12\le m \le \frac12 +4\ep_0\label{asH5}\\
&\|\nab^{\frac 12+m}(\pr^2 H)\|_{L^\infty_{[0,t_*]} L_x^2}\les 
\la^{-\f12 -4\ep_0}\quad
\mbox{for}\quad  -\frac12 + 4\ep_0 \le m \label{asH5'}\\
&\|\nab^m \big({H^{\a\b}}\pr_\a\pr_\b H\big)\|_{L^1_{[0,t_*]} L_x^\infty}\les 
\la^{-1-8\eps_0},
\label{asH4}\\
&\|\nab^m \big(\nab^{\f12}\ric(H)\big)\|_{L^\infty_{[0,t_*]} L_x^2}\les 
\la^{-1},
\label{asH800}\\
&\|\nab^m \ric(H)\|_{L^1_{[0,t_*]} L_x^\infty}\les 
\la^{-1-8\ep_0}.
\label{asH6}
\end{align} 
\end{proposition}

\section{ Generalized energy estimates
and  the Boundedness theorem}

Consider the    Lorentz metric $H=H_{(\la)}$ as in 
\eqref{Hla}  verifying, in particular, 
the properties of proposition \ref{propH} in the region
$[0,t_*]\times \rrrr^3$, $t_{*}\le \la^{1-8\ep_0}$. We denote 
by $D$ the compatible  covariant derivative and by $\nab$ the induced
covariant differentiation on $\Si_t$. We denote by $T$ the 
future oriented unit normal to $\Si_t$ and by $k$ the second fundamental
form.

  Associated 
to $H$ we have the energy momentum tensor of $\square_{H}$,
\be{enmom}
 Q_{\mu\nu}=Q[\psi]_{\mu\nu}=\pr_\mu\psi\pr_\nu\psi-
\f12 H_{\mu\nu}(H^{\a\b}\pr_\a\psi\pr_\b\psi).
\end{equation}
The energy density associated to an arbitrary timelike  vectorfield $K$
is given by $Q(K, T)$.
We consider also the modified  energy density,
\be{moden}
\Qb(K,T)=\Qb[\psi](K,T)=Q[\psi](K,T)+2t
\psi T(\psi)-\psi^2T(t).
\end{equation}
and the total conformal energy,
\be{confen}
\QQ[\psi](t)=\int_{\Si_t}\Qb[\psi](K,T).
\end{equation}
We recall below the statement of the main  generalized energy estimate
 upon which we rely.
\begin{proposition} Let $K$ be an arbitrary vectorfield 
 with deformation tensor
$$\pik_{\mu\nu}=\Lie_KH_{\mu\nu}=D_{\mu} K_{\nu}+D_{\nu} K_{\mu}$$
and $\psi$ a solution of $\square_H\psi =0$. Then
\be{genenergy}
\QQ[\psi](t)=\QQ[\psi](t_0)-\f12 \int_{[t_0,t]\times \rrrr^3}Q^{\a\b} \pikt_{\a\b}
+\frac{1}{4} \int_{[t_0,t]\times \rrrr^3} \psi^2\square_H\Omega
 \end{equation}
where
\be{pikt}
\pikt=\pik-\Omega H
\end{equation}
and $\Omega$ an arbitrary function.
\label{confenergyP}
\end{proposition}
\begin{remark} In the particular case of the Minkowski
spacetime we can choose $K$ to be the conformal timelike
 Killing vectorfield $$K=\f12 \bigg((t+r)^2(\pr_t+\pr_r)+(t-r)^2(\pr_t-\pr_r)\bigg).$$
In his case we can choose $\Omega =4t$ and obtain the total conservation law,
$$\QQ[\psi](t)=\QQ[\psi](t_0).$$
This conservation law can be used to get the desired decay estimate
for the free wave equation, see \cite{Kl}.
\end{remark}
As in  \cite{Kl-Rod}  we construct a special  vectorfield $K$
whose modified deformation tensor $\pikt$ is such that
we can control
the error terms 
$$\int_{[t_0,t]\times \rrrr^3}Q^{\a\b} \pikt_{\a\b}
+\frac{1}{4} \int_{[t_0,t]\times \rrrr^3} \psi^2\square_H\Omega.$$ 
As in \cite{Kl-Rod} we set\footnote{Observe that this definition of $K$ 
differs from the one in \cite{Kl-Rod} by an important factor of $n$.}
\be{Kvect}
 K=\f12 n(u^2\Lb+\ub^2L)
\end{equation}
with $u,\ub, L,\Lb$ defined as follows:

\begin{itemize}
\item \textit{  Optical function $u$}

 This is 
an outgoing  solution of the Eikonal equation
\be{it1}
H^{\a\b}\pr_\a u \pr_\b u=0
\end{equation}
 with initial conditions
$u(\Ga_t)=t$ on the time axis.  The time axis is
 defined as the integral curve of the forward  unit normal 
$T$ to the hypersurfaces $\Si_t$.  The point $\Ga_t$ is
the intersection between $\Ga$ and $\Si_t$. The level
surfaces of u, denoted $C_u$ are outgoing  null
 cones with vertices
on the time axis. Clearly,
\be{it1'}
T(u)=|\nab u|_h
\end{equation}
where $h$ is metric induced by $H$ on $\Si_t$, 
$|\nab u|_h^2=\sum_{i=1}^3|e_i(u)|^2$ relative to an orthonormal
frame $e_i$ on $\Si_t$.
\item  \textit{Canonical null pair $L,\Lb$}
 \be{it2}
L=bL'=T+N, \qquad \Lb=2T-L=T-N
\end{equation}
with $L'=-H^{\a\b}\pr_\b u \pr_\a$ the geodesic null generator 
 of $C_u$,   $b$ the lapse of the null foliation(or shortly null lapse)
defined by
\be{it3}
b^{-1}=-<L', T>=T(u),
\end{equation} 
and $N$ exterior unit normal, along $\Si_t$, to the surfaces
 $S_{t,u}$,  i.e. the surfaces of intersection
between $\Si_t$ and  $C_u$. We shall also
use the notation $$e_3=\Lb,\qquad  e_4=L$$

\item  \textit{The function  $\ub=-u+2t$.}
\item
\textit{The $S_{t,u}$ foliation}

The  intersection between the level hypersurfaces\footnote{The level hypersurfaces
of $u$ are outgoing  null cones $C_u$  with vertices on the time axis $\Ga_t$.}
and $u$  form compact 2- Riemannian surfaces denoted by $S_{t,u}$. We define 
$r(t,u)$ by the formula 
 Area($S_{t,u}$)$=4\pi r^2$. 
 We denote by
$\nabb$ the induced covariant derivative on $S_{t,u}$.   A vectorfield $X$
is called $S$-tangent if it is tangent to $S_{t,u}$ at every point.
 Given an $S$-tangent vectorfield  $X$  we denote by $\nabb_N X$ the projection
on $S_{t,u}$ of $\nab_N X$.   
\end{itemize}

With the help of these constructions
the  proof of the $L^2-L^\infty$ decay estimate stated in 
 theorem \ref{L2decay} can be
reduced to the following:

\begin{theorem}[Boundedness Theorem]
Consider the    Lorentz metric $H=H_{(\la)}$ as in 
\eqref{Hla}  verifying, in particular, 
the properties of proposition \ref{propH} in the region
$[0,t_*]\times \rrrr^3$, $t_{*}\le \la^{1-8\ep_0}$.
Let $\psi$ be a solution of the wave equation 
\be{psigeom}
\square_H \psi = \frac {1}{\sqrt{|H|}} 
\pr_\a \big(H^{\a\b} \sqrt{|H|}\pr_\b\psi \big)=0 
\end{equation}
 with
initial data  $\psi[t_0]$, at $t=t_0>2$,     supported in 
the geodesic ball $B_{{\half}} (0)$. Let $\DD_{u'}$ be the region
determined by $\,\,u>u'$ in the slab $[0,t_*]\times
\rr^3$.
 For all  $\,\, t_0\le t\le t_*$, $\psi(t)$ is supported 
in $\DD_{t_0-1}\subset\DD_0$ and 
$$
 \QQ[\psi](t)\les  \QQ[\psi](t_0).
$$
\label{BThrm}
\end{theorem}
We consider also  the auxiliary energy type  quantity,
\be{it7}
\EE[\psi](t)=\EEi[\psi](t)+\EEe[\psi](t)
\end{equation} where,
\beaa
\EEi[\psi](t)&=&\int_{\Si_t}(1-\ze)(t^2|\pr\psi|^2+\psi^2)\\
\EEe[\psi](t)&=&\int_{\Si_t}\,\zeta \,( \ub^2 (L\psi )^2 +
u^2 (\Lb\psi )^2 + \ub^2 |\nabb \psi |^2 + \psi^2 ).
\eeaa
with $\ze$ is a smooth cut-off function equal to $1$ in the wave
zone region $u\le \frac{t}{2}$.

In the proof of  theorem \ref{BThrm} we need  the following comparison
between the quantity $\QQ(t)$ and the auxiliary norm
$\EE(t)=\EE[\psi](t)$. 

\begin{theorem}[Comparison Theorem]
Under the same assumptions as in theorem \ref{BThrm}
we have,  for any  $1\le t\le t_*$,
$$
\EE[\psi](t)\les \QQ[\psi](t).
$$
\label{EThrm}
\end{theorem}
\section{Asymptotics Theorem and other geometric tools}
In this section we record the crucial properties of all
the important  geometric objects
associated  to our spacetime foliations $\Si_t$, $C_u$ and 
$S_{t,u}$ introduced above. Most of the  results of this section
will be proved only in the second part of this work. 

We start with some simple facts concerning the parameters
of the foliation $\Si_t$ relative to the spacetime geometry
associated to the metric $H=H_\la$. 

\medn
\textit{The $\Si_t$ foliation}\quad 
 Recall, see \eqref{decompHH},  that  the parameters of the $\Si_t$
foliation are given by $n, v$, the induced  metric $h$ and the second
fundamental form $k_{ij}$, according to the decomposition,
\be{2.3}
H=-n^2 dt^2+h_{ij}(dx^i+v^i dt)\otimes(dx^j+v^j dt),
 \end{equation} 
with $h_{ij}$ the induced Riemannian  metric on $\Si_t$,
$n$ the lapse and $v=v^i\pr_i$ the shift of $H$. Denoting by
$T$ the unit, future oriented, normal to $\Si_t$ and $k$
the second fundamental form $k_{ij}=-<\dd_iT, \pr_j>$ we 
find,
\begin{align}
 &\pr_t=n T+v,\qquad <\pr_t, v>=0\nn \\
 &k_{ij}= -\frac{1}{2}\Lie_T H_{\,ij}=-12 n^{-1}(\pr_t
 h_{ij}-\Lie_v h_{\,ij})\label{lit}
\end{align}
with $\Lie_X$ denoting the Lie derivative with respect to the 
vectorfield $X$.
 We also have 
the following, see \eqref{timelike}, \eqref{psithH}, and \eqref{psith} in section 8:
\be{4.1}
c|\xi|^2\le h_{ij} \xi^i\xi^j\le c^{-1}|\xi|^2, \qquad c \le n^2-|v|_h^2
\end{equation}
for some $c>0$.
Also
\bea
n, |v|&\les& 1\\
|\pr n|+|\pr v|+|\pr h|+|k| &\les& |\pr H| 
\eea

\textit{$S_{t,u}$- foliation}\quad  We define the Ricci coefficients associated
to the $S_{t,u}$ foliation and null pair $L,\Lb$.
\begin{definition} 
Using an arbitrary orthonormal  frame  $(e_{A})_{A=1,2}$
on $S_{t,u}$ we
 define the following tensors on the surfaces $S_{t,u}$ 
\begin{alignat}{2}
&\chi_{AB}=<\dd_A e_4,e_B>, &\quad 
&\chib_{AB}=<\dd_A e_3,e_B>,\label{chi}\nn\\
&\eta_{A}={\half} <\dd_{\bf 3} e_4,e_A>,&\quad
&\etab_{A}={\half} <\dd_{\bf 4} e_3,e_A>,\label{eta}\\
&\xib_{A}={\half} <\dd_{\bf 3} e_3,e_A>.\nn
\end{alignat}
\label{Riccicoeficients}
\end{definition}
 Using the parameters $n,v, k$ of the $\Si_t$ foliation we find(see 
\cite{Einst2} and \cite{Kl-Rod}),
\beaa
\chib_{AB}&=&-\chi_{AB}-2k_{AB}\\
\etab_A&=&-k_{AN}+ n^{-1}\nabb_A n\\
\xib_A&=&k_{AN}-\eta_A+ n^{-1}\nabb_A n\\
\eta_A&=&b^{-1}\nabb_A b+k_{AN}
\eeaa
Thus all the Ricci coefficients can be 
expressed in terms of $k_{ij}$, $n$,  the scalar function $b$ and, most
important, the Ricci coefficients $\chi$ and $\eta$.

We shall also  denote
by $\th_{AB}=<\nabb_A N, e_B> $ the second fundamental form of  $S_{t,u}$
relative to $\Si_t$.  It is easy to check that
$$\chi_{AB}=-k_{AB}+\th_{AB}.$$

We consider the parameters $b$, $\trch$, $\chih$ and $\eta$
associated to the $S_{t,u}$ foliation according to \eqref{it3} and \eqref{eta}.
For convenience we shall introduce the quantity:
\be{Theta}
\Theta=|\trch-\frac 2r|+ |\trch -\frac 2{n(t-u)}|+  |\chih|+|\eta|
\end{equation}
\begin{remark} Strictly speaking we need only one of the two quantities
$|\trch-\frac 2r|,\, |\trch -\frac 2{n(t-u)}|$ in the expression above.
Indeed we show in \cite{Einst2}  that these two are comparable.
\end{remark}
\begin{remark}
Simple calculations based on the definition \ref{Riccicoeficients},
see also Ricci equations in section 2 of \cite{Einst2},
 allow us to derive the following:
\be{Dframe}
|DL|, |D\Lb|, |\nab N|\les r^{-1}+\Theta +|\pr H|
\end{equation}
\end{remark}
\begin{remark} We shall make use of the following simple
commutation estimates, see lemma 3.5 in \cite{Einst2}, 
\be{comm}
|(\nabb_N\nabb-\nabb\nabb_N)f|\les \big(r^{-1}+\Theta +|\pr H|\big)|\nab f|
\end{equation}
\end{remark}
We state below the crucial theorem which establishes the desired
asymptotic behavior of these quantities  relative to $\la$.
\begin{theorem}[Asymptotics Theorem]
In  the spacetime region ${\cal D}_0$( see theorem \ref{BThrm})  the quantities $ b$,
$\Theta$
satisfy the following estimates:

\bea
|b-n|&\les &\la^{-4\ep_0}\\
\|\Theta\|_{L^2_t L^\infty_x}&\les& \la^{-\frac 12-3\eps_0},\label{trih1}\\
\|\Theta\|_{L^q(\stu)} &\les &\la^{-3\eps_0}.
\label{trih2}
\eea
 
In addition, in the exterior region $u\le t/2$, 
\begin{equation}
\begin{split}
&\|\Theta\|_{L^\infty(\stu)}\les t^{-1} \la^{-\eps_0} + \la^{\ep}\|\pr
H(t)\|_{L^\infty_x}.
\end{split}
\label{trih3}
\end{equation}
for an arbitrarily small  $\ep>0$.

We also have the following estimates for the derivatives of $\trch$:
\begin{align}
& \|\sup_{u\le \frac t2}
\|\Lb(\trch -\frac 2{r})\|_{L^2(\stu)}\|_{L^1_t} + \|\sup_{u\le \frac t2}
\|\Lb(\trch - \frac 2{n(t-u)})\|_{L^2(\stu)}\|_{L^1_t}\le \la^{-3\eps_0 },
\label{2trih6}\\
&\|\sup_{u\le \frac t2}\|\nabb \trch \|_{L^2(\stu)}\|_{L^1_t} +
\|\sup_{u\le \frac t2}\|\nabb \big (\trch -\frac 2{n(t-u)}\big)
\|_{L^2(\stu)}\|_{L^1_t}\le \la^{-3\eps_0 }
\label{2trih8}
\end{align}
In addition we also have  weak  estimates of the form,
\be{lastz}
\sup_{u\le \frac t2}\|(\nabb, \Lb)\big( \trch -\frac
2{n(t-u)}\big)\|_{L^\infty(\stu)}\les \la^C
\end{equation}
for some large value of $C$.

We also have the following comparison between the functions $r$ and
$t-u$,
\be{rtu}
c^{-1}\le \frac{r}{t-u}\le c
\end{equation}
\label{Asymptotics}
\end{theorem}

The proof of the Asymptotics Theorem is 
truly at  the heart of this work and it is quite involved.
 Our second  paper \cite{Einst2} is almost  entirely dedicated to it.
\begin{remark} Observe that the estimate \eqref{trih1} holds 
true also for $\pr H$. We shall  show, see \cite{Einst2} proposition 7.4,
that the $\pr H$ also verifies the estimate \eqref{trih2}.
Thus  we can incorporate the term $|\pr H|$ in the definition
\eqref{Theta} of $\Theta$.
\be{Theta'}
\Theta=|\trch-\frac 2r|+ |\trch -\frac 2{n(t-u)}|+  |\chih|+|\eta|+|\pr H|
\end{equation}
 We shall do this freely throughout this paper.
\end{remark}

 The proof of the next
proposition will be  delayed to \cite{Einst3}, see also \cite{Kl-Rod}.

\begin{proposition} Let $S_{t,u}$ be a fixed surface in  $\Si_t\cap {\cal D}_0$.

\medn
{\bf i.)}\,\,\,{\bf Isoperimetric inequality}\quad For any smooth function $f: S_{t,u}\to
\Bbb R$ we have  the following isoperimetric 
inequality:
\begin{equation}
\Bigl (\il_{S_{t,u}} |f|^2\,\Bigr )^{{\half} }
\les \il_{S_{t,u}} (|\nabb f| + |\tr\theta | |f|).     
\label{asob3}
\end{equation}
\medn
{\bf ii.)}\,\,\,{\bf Sobolev Inequality}\quad    For any $\de\in (0,1)$ and
$p$ from the interval $p\in (2,\infty]$

\begin{equation}
\begin{split}
\sup_{S_{t,u}} |f|\,\les&\, r^{\frac {\eps(p-2)}{2p+\de(p-2)}}
\bigl (\il_{S_{t,u}} (|\nabb f|^2 + r^{-2} |f|^2)
 \bigr)^{\frac 12-\frac{\de p}{2p+\de(p-2)}}\\
&\bigg[{\il_{S_{t,u}} (|\nabb f|^p + r^{-p} |f|^p)}
\bigg]^{\frac{2\de}{2p+\de(p-2)}}, 
\end{split}    
\label{asob4}
\end{equation}
\medn {\bf iii.)}\,\,\,{\bf Trace Inequality}\quad 
For  an arbitrary function $f:\Si_t\rightarrow \rr$ such that 
$f\in H^{\f12+\ep}(\rr^3)$ we have, 
\be{traceineq}
\|f\|_{L^2(S_{t,u})}\le\|\pr^{\f12 +\ep}f\|_{L^2(\Si_t)}+\|\pr^{\f12
-\ep}f\|_{L^2(\Si_t)}.
\end{equation}
More generally, for any $q\in [2,\infty)$
\be{traceineq2}
\|f\|_{L^q(S_{t,u})}\le\|\pr^{\frac 32 -\frac 2q +\ep}f\|_{L^2(\Si_t)}+\|\pr^{\frac 
32 -\frac 2q -\ep}f\|_{L^2(\Si_t)}.
\end{equation}
Also, considering the region $\mbox{Ext\,}_t=\Si_t\cap\{0\le u\le \frac t2\}$,
 we have the following:
\be{traceineq22}
\|f\|_{L^2(S_{t,u})}^2\le\|N(f)\|_{L^2(\mbox{Ext\,}_t)}\|f\|_{L^2(\mbox{Ext\,}_t)}+
\frac{1}{t}
\|f\|_{L^2(\mbox{Ext\,}_t)}.
\end{equation}
\label{Triso}
\end{proposition} 
We shall make use of  the following, see  lemma 6.3 in \cite{Kl-Rod}.
\begin{proposition}
\label{Polar}
The following inequality holds for all $t\in [1,t_*]$ and $2<p<\infty$:
\begin{equation}
\il_{\Si_t} V^2 w^2 \le 
t^{\frac{2}{p}} \sup_{ u } \|V\|^2_{L^{2p'}(S_{t,u})}
\int_{\Si_t}\big( |\nabb w|^2 +r^{-2} |w|^2\big).
\label{polar100}
\end{equation}
where $p'$ is the exponent dual to $p$.

We shall also make use of the form,
\be{polar101}
\il_{\Si_t} V^2 w^2 \le 
t^{\frac{2}{p}} \|V\|_{L_x^\infty}^{\frac{2}{p}} \sup_{ u }
\|V\|^{\frac{2}{p'}}_{L^{2}(S_{t,u})}
\int_{\Si_t}\big( |\nabb w|^2 +r^{-2} |w|^2\big).
\end{equation}
In particular, if $\|V\|_{L_x^\infty}$ is bounded by
some positive  power of $\la$, and  we restrict ourselves
to the exterior region $\Ext$, we deduce that for 
every $\varep>0$ and some constant  $C$
\begin{equation}
\il_{Ext_{t}} V^2 w^2 \le t^{-2}
\la^{C\varep} \sup_{0\le u\le t/2 } \|V\|^{2-\varep}_{L^2(S_{t,u})} \EE [w](t).
\label{polar}
\end{equation}
\end{proposition}

\begin{proof} The proof is straightforward and relies only on the
isoperimetric inequality \eqref{asob3}, see also 6.1. in  \cite{Kl-Rod}.
\end{proof}

\section{Proof of the Boundedness Theorem}
We first calculate the components of the modified\footnote{ corresponding
to the choice $\Om =4t$.}  deformation
tensor
$\pit=\pikt=\pik-4t H$
of  our vectorfield $K= \frac 12 n(\ub^2 L + u^2 \Lb)$.
Recall that $\ub=2t-u$ and $\Lb=-L+2T$, thus 
\beaa
\Lb\,(u^2)&=&4u\, b^{-1},\\
L(\,\ub^2)&=&4\ub\, n^{-1},\\
\Lb(\,\ub^2)&=&4\ub\, (n^{-1}-b^{-1}).
\eeaa

Proceeding as in section 6.1 of \cite{Kl-Rod} we  calculate
 the null components of $\pit=\pikt$ relative\footnote{We say that $(e_{i})_{1=1,2,3,4}$
 forms a null frame.}
to 
 $e_4=L, e_3=\Lb$ and $(e_A)_{A=1,2}$ an arbitrary  orthonormal
frame on $S_{t,u}$ find,
\begin{align}
&\pit_{\bf 44} = 2u^2 \,n(\bk_{NN}- n^{-1} e_4(n)),\label{dek1}\nn\\
&\pit_{\bf 34}= 4u n(n^{-1}-b^{-1})+ \ub^2 n\big (\bk_{NN}- n^{-1}e_4 (n)\big ) + 
u^2 n \big(\bk_{NN}-n^{-1} e_3(n)\big),\label{dek2}\nn\\
&\pit_{\bf 33}= -8\ub\, n(n^{-1}-b^{-1})-2\ub^2\, n
\big(\bk_{NN} + n^{-1} e_3(n)\big ),\label{dek3}\nn\\
&\pit_{{\bf 3}A}= \ub^2\, n(\eta_A + k_{AN} - n^{-1} \nabb_A n) + u^2 n \xib_A ,
\label{dek4}\\
&\pit_{{\bf 4}A}= u^2 \,n (\etab_A-k_{AN}- n^{-1} \nabb_A n),\label{dek5}\nn\\
&\pit_{AB}= 2tn(t-u)\big (\trch - \frac 2{n(t-u)}\big)\de_{AB} + 
4tn(t-u)\chih_{AB}- 2u^2 n k_{AB}\label{dek6}\nn
\end{align}

The  following proposition concerning 
 the behavior of the null  components of $\pit$
is an immediate consequence of  the above formulae and 
the Asymptotics Theorem stated above.
\begin{proposition}
\label{deftensestim}
\begin{alignat}{2}
&\|u^{-2} \pit_{44}\|_{L^1_t L^\infty_x} \les \la^{-3\eps_0}, &\quad &\quad
\|(u\ub)^{-1} \pit_{34}\|_{L^1_t L^\infty_x} \les \la^{-3\eps_0},\nn \\
&\|(\ub)^{-2} \pit_{33}\|_{L^1_t L^\infty_x} \les \la^{-3\eps_0},&\quad &\quad
\|(\ub)^{-2} \pit_{3A}\|_{L^1_t L^\infty_x} \les \la^{-3\eps_0},\nn \\
&\|(u)^{-2} \pit_{4A}\|_{L^1_t L^\infty_x} \les \la^{-3\eps_0},&\quad &\quad
\|(\ub)^{-2} \pit_{AB}\|_{L^1_t L^\infty_x} \les \la^{-3\eps_0}.\nn
\end{alignat}

\end{proposition}
The  proof of  the Boundedness theorem
relies on the generalized energy identity \eqref{genenergy}
with $K={\half} n\big( u^2\Lb+\ub^2 L\big)$ and $\Omega=4t$.
Thus,
\bea
\QQ[\psi](t)&=&\QQ[\psi](t_0)-\f12 
\int_{[t_0,t]\times \rrrr^3}Q^{\a\b}
\pikt_{\a\b} + \int_{[t_0,t]\times \rrrr^3} \psi^2\square_H t\nn\\
&=&\QQ[\psi](t_0)-\f12{\cal J}+{\cal Y}
\label{JY}
\eea
Observe that we can decompose:
$$
\aligned
\ & {\cal J}=\il_{[t_0,t]\times\rr^3}  Q^{\a\b}[\psi] \pit_{\a\b} = 
\il_{[t_0,t]\times\rr^3}\bigg(\frac 14 \pit_{33}(L\psi)^2 + 
\frac 14\pit_{44} (\Lb \psi)^2 +
{\half} \pit_{34} |\nabb \psi |^2 \\  \ & -\pit_{{4}A}\Lb \psi \nabb_A \psi -
\pit_{{3}A}L \psi \nabb_A \psi + \pit_{AB} \nabb_A \psi \nabb_B \psi + 
\tr \pit ({\half}\Lb\psi L\psi -|\nabb\psi|^2)\bigg).
\endaligned
$$
 Consider, for example, $I=\il_{[t_0,t]\times\rr^3} \pit_{{4}A}\Lb \psi \nabb_A
\psi$. We can estimate it  as follows :
$$
\aligned
I
 &\le {\half} \il_{[t_0,t]\times\rr^3} |(\ub u)^{-1}\pit_{{4}A}|\bigg( u^2 (\Lb \psi)^2 +
\ub^2 (\nabb_A \psi)^2 \bigg)\\ &\le \il_{t_0}^t  \|(\ub u)^{-1}\pit_{{4}A}\|_{L^\infty_x}
\EE[\psi](\tau)\,d\tau
\endaligned
$$
Making  use of the comparison theorem 
and the estimate $ \|(\ub u)^{-1}\pit_{{4}A}\|_{L^1_t L^\infty_x}\les \la^{-3\eps_0}$ we
infer that,
$$I\le \il_{t_0}^t  \|(\ub u)^{-1}\pit_{{4}A}\|_{L^\infty_x}
\QQ[\psi](\tau)\,d\tau\les\la^{-3\eps_0}\sup_{[t_0,t]} \QQ[\psi](\tau)
$$ 
 We can proceed in the same manner with all the terms  of ${\cal J}$ with the exception
of $\il_{[t_0,t]\times\rr^3}\tr \pit\, \Lb\psi\,  L\psi$. Observe 
that\footnote{We use $tr$ here to denote the trace relative to the 
surfaces $S_{t,u}$. Thus $\tr k=\de^{AB}k_{AB}$. We use  $\mbox{Tr}k=h^{ij}k_{ij}$
 to denote the usual trace of $k$ with respect to $\Si_t$.}
$$
\tr \pit =\de^{AB} \pit_{AB} = 2tn(t-u)\big(\trch -\frac 2{n(t-u)}\big)
- 2 u^2 n \tr k
$$
$$
\il_{[t_0,t]\times\rr^3}|u^2 n \tr k L\psi \Lb \psi |\le 
{\half} \il_{[t_0,t]\times\rr^3} |\tr k| \big (u^2 (L\psi)^2 + u^2 (\Lb \psi)^2\big)
\les \il_{t_0}^t \|\pr H\|_{L^\infty_x} \EE[\psi](\tau)\,d\tau 
$$
Since $\|\pr H\|_{L^1_t L^\infty_x} \les \la^{-4\eps_0}$, this term can be 
treated in the same manner as $ I$. We are thus left  with the integral
$$
{\cal B}=\il_{[t_0,t]\times\rr^3} 2tn(t-u)\big(\trch -\frac 2{n(t-u)}\big)\Lb\psi L\psi
$$
All other terms ${\cal J}-{\cal B}$ can be estimated in precisely
the same  manner, using  the comparison theorem and the estimates of theorem
\ref{deftensestim},  by
\be{good}
{\cal J}-{\cal B}\lesssim \la^{-3\ep_0} \sup_{[t_0,t]} \QQ[\psi](\tau)
\end{equation}
To estimate the remaining term ${\cal B}$ requires a more involved 
argument. In fact we shall need  more information concerning 
the geometry of the null cones $C_u$ and surfaces $S_{t,u}$.

Denote $Ext_t$ the exterior region $Ext_t = \{0\le u \le t/2\}$.
Let $\zeta$ be a smooth cut-off function with support in $Ext_t$.
Observe that 
\begin{equation}
\il_{\Si_t} 
\big (t^2 (\pr \psi)^2 +\psi^2\big)(1-\zeta) \les \il_{\Si_t} (1-\zeta) \Qb[\psi](t)
\label{polar1}
\end{equation}
We can split the remaining integral 
\beaa
{\cal B}&=&{\cal B}^{i}+{\cal B}^{e}
\\{\cal B}^{i} &=& \il_{[t_0,t]\times\rr^3} 2tn(t-u)\big(\trch -\frac 2{n(t-u)}\big)
L\psi
\Lb \psi\, (1-\zeta)\\
 {\cal B}^{e}&=&\il_{[t_0,t]\times\rr^3} 2tn(t-u)\big(\trch -\frac 2{n(t-u)}\big) L\psi
\Lb
\psi\,\zeta
\eeaa
With the help of \eqref{polar1} the first integral can be estimated as follows:
\beaa
{\cal B}^{i}&\les& \il_{[t_0,t]\times\rr^3} |\trch -\frac 2{n(\tau-u)}|\, \,\tau^2
\,(\pr\psi)^2 (1-\zeta)\\ 
&\les &\il_{t_0}^t \|\trch -\frac 2{n(\tau-u)}\|_{L^\infty_x} \,\,
\Qb[\psi](\tau)\,d\tau
\eeaa
In view of the estimate 
$\|\trch -\frac 2{n(t-u)}\|_{L^1_t L^\infty_x} \les \la^{-3\eps_0}$,
given by the Asymptotics Theorem \eqref{Asymptotics} we infer that,
$${\cal B}^{i}\lesssim \la^{-3\ep_0} \sup_{[t_0,t]} \QQ[\psi](\tau)$$

Therefore, it remains to estimate ${\cal B}^e$.

 According to the Asymptotics Theorem 
the quantity $z= \trch -\frac 2{n(t-u)}$ verifies the following estimates:
\begin{alignat}{2}
&\|z\|_{L^2_t L^\infty_x} \les \la^{-\f12-3\eps_0},
&\qquad &\|z\|_{L^2(\stu)}\les \la^{-2\eps_0},\label{polar2}\\
&\|\sup_{u\le \frac t2}\|\nabb z\|
_{L^2(\stu)}\|_{L^2_t} \les \la^{-\f12 -3\eps_0}, &\qquad & 
\|\sup_{u\le \frac t2}\|\Lb\ z\|
_{L^2(\stu)}\|_{L^2_t} \les \la^{- \f12 -3\eps_0}.\label{polar3}
\end{alignat} 
\begin{remark} The same estimates hold true
if we replace $\trch -\frac 2{n(t-u)}$ by $\trch -\frac 2r$.
\label{rempolar}
\end{remark}
It would therefore suffice to prove the following result.
Using the estimates  \eqref{polar2}--\eqref{polar3} we shall
prove that:
\be{badt}
{\cal B}^e=\il_{[t_0,t]\times\rr^3} 2tn(t-u)z L\psi
\Lb
\psi\,\zeta \les \la^{-\ep_0}\sup_{[t_0,t]} \QQ[\psi](\tau)
\end{equation}
To prove \eqref{badt} we need to rely on the fact
that $\psi$ is  a solution of the wave equation $\square_H\psi=0$.
We  shall also make use of  the following   
 standard integration by parts formulae\footnote{These are 
simple adaptations of the formulae in   lemma 6.2.,
\cite{Kl-Rod}. }, 
 \begin{equation}
\il_{\Si_t} F N(G) = - \il_{\Si_t} \bigg(N(F) + \big(\tr\th + n^{-1}N(n)\big) F\bigg)G,
\label{ipar3}
\end{equation} 
where  $N$ is  the unit normal to
$\stu$.  

If  $Y$ is a vectorfield in $T \Si_t$ tangent to 
$S_{t,u}$ then
\begin{equation}
\il_{\Si_t} F \divv Y = - \il_{\Si_t} \bigg(\nabb F + (b^{-1}\nabb b + n^{-1}\nabb n)
F\bigg)\cdot Y.
\label{ipar4}
\end{equation}
It is also not difficult to verify that
\begin{equation}
\il_{[t_0,t]\times \rr^3} F T(G) = -\il_{[t_0,t]\times \rr^3} 
(T(F)  + Tr k +\div v)G +
\il_{\Si_t} F G - \il_{\Si_{t_0}} FG 
\label{ipar5}
\end{equation}
Writing   $\Lb=T-N$  we   integrate by parts and express   the integral 
${\cal B}^e$
in the form,
\bea
{\cal B}^e & = &-I_1 +I_2 +I_3 -I_4\\
I_1&=&\il_{[t_0,t]\times\rr^3}\zeta\,  nt(t-u) z \,
(\Lb L\psi)\, \psi \,\nn\\
I_2&=&\il_{[t_0,t]\times\rr^3}\Bigl (-\Lb (\zeta\, nt(t-u) z) +
(\tr\th +n^{-1}N(n)- \Tr k-\div v )\zeta\, nt(t-u) z\Bigr ) L\psi\, \psi\nn\\
I_3&=&\il_{\Si_t}\zeta\, nt(t-u) z \,L\psi \psi\nn\\
I_4&=&\il_{\Si_{t_0}}\zeta\, nt(t-u) z \,L\psi \psi\nn
\eea

We first handle the boundary terms $I_3$, $I_4$.
With the help of proposition \ref{Polar}( which we can apply 
 in view of the estimates \eqref{trih3} for $\Theta$ as well
as the estimate \eqref{asH2} for $\pr H$.)
we have
$$ \|n(t-u) z \psi\|_{L^2(Ext_t)} \les  \la^{C\eps} \sup_{u\le \frac t2} \|n
z\|_{L^2(\stu)}^{1-\ep/2}\EE^{{\half}}[\psi](t).  $$
Therefore,
$$
\aligned
\il_{\Si_t}| \zeta\, nt(t-u) z \,L\psi \psi |& \les 
\il_{Ext_t}|n(t-u) z \,t L\psi \psi |\\
 &\les\|t\, L\psi\|_{L^2(\Si_t)}\|n(t-u)z\psi\|_{L^2(Ext_t)}\\
 &\les 
\|n(t-u)\, z \psi\|_{L^2(Ext_t)} \EE^{{\half}}[\psi](t)\\
& \les \la^{C\eps} \sup_{s\ge \frac t2} \|\,n z\|_{L^2(\stu)}^{1-\ep/2}
\EE[\psi](t)\les \la^{-\eps_0} \EE[\psi](t).
\endaligned
$$
The last inequality followed from the boundness of $n$ and
\eqref{polar2}. Similar estimate holds for the second boundary
term $I_4$.

To estimate $I_2$  we observe that, as an immediate consequence
of theorem \ref{Asymptotics}, we have
$$|\Lb(t)|,\,\, |\Lb(t-u)|\les 1,\qquad |\Lb (\zeta)|\les t^{-1}$$

Denoting  
$$\Theta(t,x)= |\trch -\frac{2}{n(t-u)}| +|\chih| +|\eta| +|\pr H|$$
we easily find,
$$
|I_2 |\les \il_{t_0}^t\il_{Ext_\tau}\bigg (\tau^2 |\Lb(z)| + \tau|z| + \tau^2 \Theta
|z|\bigg)   | L\psi\, \psi |\,d\tau
$$

To treat the term involving $\Lb(z)$
we proceed as in the case of $I_1$;  We estimate
 $\il_{Ext_\tau}\tau^2 |\Lb(z)|| L\psi\, \psi |\,d\tau$
by Cauchy-Schwartz followed by an application of proposition \ref{Polar}.
 The space integral of the  other two terms  can be estimated
as follows:
$$\il_{Ext_\tau}\big ( \tau|z| + \tau^2 \Theta
|z|\big)   | L\psi\, \psi |\,d\tau\le \big(\|z\|_{L_x^\infty}+
\tau\|\Theta\|_{L_x^\infty}\|z\|_{L_x^\infty}\big) \EE[\psi](\tau).
$$
Consequently, using the inequalities
  \eqref{polar2}-\eqref{polar3} for $z$ ( as well as  the weak estimate \eqref{lastz})
and the estimates for $\Theta$ from the Asymptotics Theorem \ref{Asymptotics}
\beaa
I_2&\les& \il_{t_0}^t \bigg( \la^{C\eps} \sup_{u\le \frac {\tau}2} 
\|\Lb(z)\|_{L^2(\stu)}^{1-\ep/2} + \|z\|_{L^\infty(\Si_\tau)} + \tau 
\|\Theta\|_{L^\infty(\Si_\tau)} 
\|z\|_{L^\infty(\Si_\tau)}\bigg ) \EE[\psi](\tau)\,d\tau\\
&\les& \la^{C'\eps} \bigg(\|\sup_{u\le \frac {\tau}2} 
\|\Lb(z)\|_{L^2(\stu)}\|_{L^1_t}^{1-\ep/2} +  \|z\|_{L^1_t L^\infty_x} + 
\la  \|\Theta\|_{L^2_t L^\infty_x} 
\|z\|_{L^2_t L^\infty_x}\bigg)\sup_{[t_0,t]}\QQ[\psi](\tau)\\
 &\les&\la^{-\eps_0}\sup_{[t_0,t]}\QQ[\psi](\tau)
\eeaa
as desired.

It remains therefore to consider $I_1$. 
We shall make  use of the fact 
that $\psi$ is a solution of the wave equation. This  allows us
to express the $ \Lb L(\psi)$  in terms of
the angular laplacian\footnote{the Laplace-Beltrami
operator on $S_{t,u}$.} $\lapp$  and lower order terms. 
Expressed relative to a  null frame the wave operator  $\square_H \psi $  takes the form
$$
\square_H \psi = H^{\a\b} \psi_{;\a\b} = -\psi_{;\bf 43} + \psi_{;AA},
$$
where $\psi_{;e_i e_j} = e_j (e_i (\psi )) - \dd_{e_i} e_j (\psi)$.
We use the Ricci formulas: 
$
\dd_{\bf 3} e_4 = 2 \eta_A e_A + \bk_{NN} e_4,$ and  
$\dd_B e_A = \nabb_B e_A + {\half} \chi_{AB} e_3 + {\half} \chib_{AB} e_4\,\,$
to derive 
\begin{equation}
\square_H \psi = - \Lb L\psi + \lapp \psi + 
2 \eta_A \nabb_A \psi + {\half}\trch \Lb \psi + 
({\half} \trchb + \bk_{NN}) L \psi .
\label{wframe}
\end{equation}
As a result of this calculation
\bea
I_1&=&\il_{[t_0,t]\times\rr^3}\zeta\, n\tau(\tau-u) z\,\Lb L\psi\,\psi  = 
\il_{[t_0,t]\times\rr^3}\zeta\, n\tau(\tau-u) z\,\lapp \psi\psi\nn\\
&+&{\half} \il_{[t_0,t]\times\rr^3}\zeta\, n\tau(\tau-u) z\,\trch (\Lb \psi)\psi\nn\\
&+& \il_{[t_0,t]\times\rr^3}\zeta\, n\tau(\tau-u) z\,\Bigl (2 \eta_A \nabb_A \psi + 
({\half} \trchb + \bk_{NN}) L \psi \Bigr )\,\psi\nn\\
& = &I_{11} + I_{12}+I_{13}. \label{ipar100}
\eea
Consider first  $I_{13}$. Taking into account that $t-u\ge \frac t2$ 

\begin{align}
|I_{13}|& \les \il_{t_0}^t\il_{Ext_\tau} \tau^2 |z| \bigg(\Theta \nabb\psi + 
(\frac 1\tau + \Theta ) L\psi \bigg )\psi \nn\\ &\les 
 \il_{t_0}^t  \bigg(\tau \|z\|_{L^\infty (\Si_\tau)} \|\Theta\|_{L^\infty(\Si_\tau)} 
+ \|z\|_{L^\infty (\Si_\tau)}\bigg)\EE[\psi](\tau)\,d\tau\nn\\
&\les \la^{-\eps_0}\sup_{[t_0,t]}\QQ[\psi](\tau)
\label{I13}
\end{align}

as before.

To estimate   $I_{12}$ we need first to integrate
once more  by parts.
$$
\aligned
I_{12} &= 
\frac 14 \il_{[t_0,t]\times\rr^3}\Bigl (-\Lb (\zeta\,n\tau(\tau-u) z\,\trch)\\
& +
(\tr\th + n^{-1} N(n)- \Tr k -\div v)\,\zeta\,n\tau(\tau-u) z \trch\Bigr ) \psi^2 \\ &+
\frac 14 \il_{\Si_t} \zeta\, n\tau(\tau-u) z\,\trch
(\psi )^2 - \frac 14 \il_{\Si_{t_0}} \zeta\,n\tau(\tau-u) z\,\trch
(\psi )^2 
\endaligned
$$
All  terms can be treated  as above. Take, for example,
the worst term involving $\Lb(\trch)$.
 Recall that 
$$
\Lb(\trch ) = \Lb (\trch -\frac 2r) + \Lb\big(\frac 2r\big)\les
 \Lb (\trch -\frac 2r) + \frac 2{r^2} + \frac 1r \Theta
$$
Thus
$$
\aligned
\il_{[t_0,t]\times\rr^3}| \zeta\,nt(t-u) z \Lb(\trch)
(\psi)^2| & \les \il_{t_0}^t\il_{Ext_\tau} \tau^2 |z|\big ( |\Lb (\trch -\frac 2r)|
+ \frac 1{\tau^2} + \frac 1\tau \Theta\big ) (\psi)^2\\ &\les 
 \il_{t_0}^t\il_{Ext_\tau} \tau^2\, |z| \,|\Lb (\trch -\frac 2r)|\, \psi^2 \\ &+ 
 \il_{t_0}^t \big (\|z\|_{L^\infty(\Si_\tau)}  + \tau  \|z\|_{L^\infty(\Si_\tau)}
\|\Theta\|_{L^\infty(\Si_\tau)}\big ) \EE[\psi](\tau) \,d\tau
\endaligned
$$
The second term has already been treated above, see \eqref{ipar100}.
To estimate the first we apply first Cauchy-Schwartz and then
 make use of  proposition
\ref{Polar},
$$
\aligned
\il_{t_0}^t\il_{Ext_\tau} \tau^2 |z| |\Lb (\trch -\frac 2r)| \psi^2&\les
\il_{t_0}^t \|\,\tau^2 \, |z|\, |\Lb (\trch -\frac 2r)| \,\psi\|_{L^2(Ext_\tau)}
\,\,\EE^{{\half}}[\psi](\tau)\,d\tau \\ &\les 
\il_{t_0}^t \la^{C\eps} \sup_{u\le \frac \tau 2} 
 \|\tau |z| |\Lb (\trch -\frac 2r)| \|_{L^2(\stu)}^{1-\ep/2}
\EE[\psi](\tau)\,d\tau
\endaligned
$$
Taking into account the  estimates  in \eqref{polar2}--\eqref{polar3}
and the Remark \ref{rempolar} we 
deduce,
\beaa
\la^{C\eps}\il_{t_0}^t  \sup_{u\le \frac \tau 2} 
 \|\tau z \Lb (\trch -\frac 2r) \|_{L^2(\stu)}^{1-\ep/2} &\les &
\la^{C'\eps}\bigg( t \|z\|_{L^2_t L^\infty_x} \| \sup_{u\le \frac t 2} 
\|\Lb (\trch -\frac 2r)| \|_{L^2(\stu)}\|_{L^2_t}\bigg)^{1-\ep/2}\\
&\les&
\la^{-\eps_0}
\eeaa
Therefore,
\be{I12}
|I_{12}|\les \la^{-\ep_0}\sup_{[t_0,t]} \QQ[\psi](\tau)
\end{equation}

Finally we  estimate $I_{11}=\il_{[t_0,t]\times\rr^3}\zeta\,nt(t-u) z \lapp\psi\,\psi$
 by integrating once more by parts as follows:
$$
\aligned
I_{11} & = -
\il_{[t_0,t]\times\rr^3} \zeta\,nt(t-u) z \,|\nabb \psi|^2 \\ &-
\il_{[t_0,t]\times\rr^3} n^{-1} b^{-1}\nabb_A (b n\,\zeta\,
nt(t-u) z ) \nabb_A\psi\,\psi .
\endaligned
$$
The first integral  on the right can be
easily  estimated 
\bea
\il_{[t_0,t]\times\rr^3} \zeta\,nt(t-u) z \,|\nabb \psi|^2\nn
&\les& \il_{t_0}^t \|z\|_{L^\infty_x} \EE[\psi](\tau)\,d\tau\nn\\
&\les&\|z\|_{L_t^1L_x^\infty}\sup_{[t_0,t]}\QQ[\psi](\tau)\nn\\
&\les& \la^{-3\ep_0}\sup_{[t_0,t]}\QQ[\psi](\tau)
\label{I111}
\eea

To estimate the second we write  schematically
$$
\nabb(b n^2\,\zeta\,t(t-u) z ) \approx t(t-u) (\nabb b ) z + 
t(t-u) \nabb z + t(t-u) z \Theta = 
t(t-u) \nabb z + t(t-u) z \Theta
$$
since $\nabb_A b = b(\eta_A - k_{AN})$.  
Thus with the help of proposition \ref{Polar}( using  also the weak estimate
\eqref{lastz}),
$$
\aligned
\il_{[t_0,t]\times\rr^3} |n^{-1} b^{-1}\nabb_A (b n^2\,\zeta\,
\tau(\tau-u) z ) \nabb_A\psi\,\psi | \les 
\il_{t_0}^t\il_{Ext_\tau} \big (\tau|\nabb z| + \tau |z||\Theta|\big )
|\tau \nabb_A \psi|\,|\psi|& \\  \les 
 \il_{t_0}^t \bigg(\la^{C\eps} \sup_{u\le  \frac \tau 2}\|\nabb z\|_{L^2(\stu)}^{1-\ep/2}  
+ \tau \|z\|_{L^\infty(\Si_\tau)} \|\Theta\|_{L^\infty(\Si_\tau)} \bigg)
\EE[\psi](\tau)\,d\tau. &
\endaligned
$$ 
Using \eqref{polar3}  once more we have,
$$
\aligned
&\il_{t_0}^t \bigg(\la^{C\eps} \sup_{u\le \frac \tau 2}\|\nabb z\|_{L^2(\stu)} ^{1-\ep/2} 
+ \tau \|z\|_{L^\infty(\Si_\tau)} \|\Theta\|_{L^\infty(\Si_\tau)} \bigg)\,d\tau
\\ &\les
\la^{C'\eps} \|\sup_{u\le \frac t2} \|\nabb z\|_{L^2(\stu)} \|_{L^1_t}^{1-\ep/2}
+ t  \|z\|_{L^2_t L^\infty_x} \|\Theta\|_{L^2_t L^\infty_x}\les \la^{-\eps_0}
\endaligned
$$
Therefore, combining with \eqref{I111} we infer that,
\be{I11}
I_{11}\les\la^{-\ep_0}\sup_{[t_0,t]}\QQ[\psi](\tau)
\end{equation}
Recalling also \eqref{I12} and  \eqref{I13}
we conclude that
\be{I1new}
I_{1}\les\la^{-\ep_0}\sup_{[t_0,t]}\QQ[\psi](\tau)
\end{equation}
Since $I_2, I_3,I_4$ and ${\cal B}^i$ have already been estimated
we finally derive,
 \be{bad}
{\cal B}\les\la^{-\ep_0}\sup_{[t_0,t]}\QQ[\psi](\tau)
\end{equation}
as desired. This combined with
\eqref{good} yields,
\be{calJ}
{\cal J}\les\la^{-\ep_0}\sup_{[t_0,t]}\QQ[\psi](\tau)
\end{equation}
Going back to the identity \eqref{JY}
we still have to estimate ${\cal Y}$.
For this we only need to observe that
$\square_Ht$ depends only on the first derivatives
of $H$. Thus also
\be{calY}
{\cal Y}\les\la^{-\ep_0}\sup_{[t_0,t]}\QQ[\psi](\tau)
\end{equation}
Therefore,
$$\sup_{[t_0,t]}\QQ[\psi](\tau)\le
\QQ[\psi](t_0)+\la^{-\ep_0}\sup_{[t_0,t]}\QQ[\psi](\tau) $$
which implies the boundedness theorem.

\section{Proof of the Comparison Theorem }
We proceed precisely as in \cite{Kl-Rod}, section 6.1. 
Define  $S$ and $\Sbb$,
\begin{equation}
S = {\half} (\ub L + u \Lb), \quad \Sbb = {\half} (\ub L - u \Lb).
\label{SSbar}
\end{equation}
Since $\ub=-u+2t$, $\Lb=T-N$, $L=T+N$
$$
\aligned
\ &t T = \frac{1}{4} (u+\ub)(L+\Lb) = S - \frac 14 (u-\ub)(L-\Lb)= S-(t-u) N,\\ 
\ & t T= {\half} t (\Lb + L) = \frac {t}{t - u} \Sbb - \frac {t^2}{t-u} N.
\endaligned
$$

Therefore, with the help of the identities \eqref{ipar3}, and $N(t)=0$,
$N(u)=-b^{-1}$
$$
\aligned
2 \il_{\Si_t}\psi t T( \psi) = &2 \il_{\Si_t}(\psi (S \psi) - {\half}
(t-u) N (\psi^2))\\ = &2 \il_{\Si_t}\psi (S \psi) +
 \il_{\Si_t}\bigg(b^{-1}+(t-u)\big(\tr\th +n^{-1} N(n)\big)\bigg)\psi^2,\\
2 \il_{\Si_t}\psi t T( \psi) = &2 \il_{\Si_t}
(\psi \frac t{t-u}(\Sbb \psi) - {\half}
\frac {t^2}{t-u} N (\psi^2))\\ =  & 2 \il_{\Si_t}\psi \frac t{t-u}(\Sbb \psi) +
 \il_{\Si_t}\frac{t^2}{(t-u)^2}\bigg(-b^{-1}+(t-u)\big(\tr\th
+n^{-1}N(n)\big)\bigg)\psi^2.
\endaligned
$$
Recall that  $\th_{AB}=\chi_{AB} + k_{AB}$.  
Recall also that  $\Theta$ was defined in \eqref{Theta'}.
$$\Theta(t,x)=|\trch -\frac{2}{r}|+ |\trch -\frac{2}{n(t-u)}|+|\chih|+|\eta|  +|\pr H|$$
Thus,
$$
\aligned
\ &2 \il_{\Si_t}\psi t T( \psi) = 2 \il_{\Si_t}(\psi (S \psi) +
 \il_{\Si_t}\bigg(b^{-1}+\frac{2}{n} +(t-u)\Theta \bigg)\psi^2,\\
\ &2 \il_{\Si_t}\psi t T( \psi) = 2 \il_{\Si_t}
(\psi \frac t{t-u}(\Sbb \psi)  +
 \il_{\Si_t}\frac{t^2}{(t-u)^2}(-b^{-1}+\frac{2}{n}+(t-u)\Theta \bigg)\psi^2.
\endaligned
$$
Recall, from the Asymptotics Theorem \ref{Asymptotics}, 
$$|b-n|\les \la^{-4\ep_0}
$$
Also, since $n$ is bounded away from zero so is $b$.
Therefore,
$$
\aligned
\ &2 \il_{\Si_t}\psi t T( \psi) = 2 \il_{\Si_t}(\psi (S \psi) +
 \il_{\Si_t}\bigg(\frac{3}{n} +(t-u)\Theta +\la^{-4\ep_0}\bigg)\psi^2,\\
\ &2 \il_{\Si_t}\psi t T( \psi) = 2 \il_{\Si_t}
(\psi \frac t{t-u}(\Sbb \psi)  +
 \il_{\Si_t}\frac{t^2}{(t-u)^2}(\frac{1}{n}+(t-u)\Theta +\la^{-4\ep_0}\bigg)\psi^2.
\endaligned
$$

Since
$$
\bar{Q}(K,T)[\psi]= \frac n4 \bigl ( \ub^2 (L\psi )^2 + 
u^2 (\Lb\psi )^2 + (\ub^2+u^2) |\nabb \psi |^2 \bigr )+ 
2t \psi T \psi - n^{-1}\psi^2,
$$
and 
$$
\frac 14 ( \ub^2 (L\psi )^2 + 
u^2 (\Lb\psi )^2 ) = {\half} \bigl ((S \psi)^2 + (\Sbb \psi )^2\bigr )
$$
we can introduce positive constants $A, B :\, A+B =2$ such that
\beaa
 \QQ[\psi](t) & = &{\half} \il_{\Si_t}
\Bigl (n(S \psi)^2 + 2 A \psi (S \psi) +  (\frac{3}{n}A-\frac 2n)\psi^2+
\big((t-u)\Theta+\la^{-4\ep_0}\big)\psi^2
\Bigr )
 \\ 
& +& {\half} \il_{\Si_t}\Bigl (n(\Sbb \psi )^2 + 2 B  \psi \frac t{t-u}(\Sbb \psi) + 
(\frac{1}{n}\frac {t^2}{(t-u)^2} B\psi^2+ \frac
{t^2}{(t-u)^2}\big((t-u)\Theta+\la^{-4\ep_0}\big)\psi^2\Bigr )\\ &+& {\half}
\il_{\Si_t}n(\ub^2+u^2) |\nabb \psi |^2 .
\eeaa
For any values of $A, B$ such that $1< A < 2$ and $0<B<1$ it is possible
to find positive constants $c_1, c_2$ such that
\beaa
n(S \psi)^2 + 2 A \psi (S \psi) +  \frac{1}{n}(3A-2)\psi^2 &\geq &c_1 \bigl 
((S \psi)^2 + \psi^2\bigr ), \\
n(\Sbb \psi)^2 + 2 B \psi \frac t{t-u}(\Sbb \psi) + \frac 1n B\frac {t^2}{(t-u)^2}
\psi^2 &\geq& c_2 \bigl ((\Sbb \psi)^2 + \frac {t^2}{(t-u)^2} \psi^2\bigr ).
\eeaa

Therefore,

\beaa
\QQ[\psi](t) &\ges &
\il_{\Si_t} ( \ub^2 (L\psi )^2 + u^2 (\Lb\psi )^2
+ (u^2+\ub^2) |\nabb \psi |^2 + (1+\frac {t^2}{(t-u)^2})\psi^2\\
& -&\il_{\Si_t}(1+ \frac {t^2}{(t-u)^2}) \bigg((t-u)\Theta+\la^{-4\ep_0}\bigg)\psi^2
\eeaa

\beaa
\QQ[\psi](t) &\ges &
\il_{\Si_t} ( \ub^2 (L\psi )^2 + u^2 (\Lb\psi )^2
+ (u^2+\ub^2) |\nabb \psi |^2 + (1+\frac {t^2}{(t-u)^2})\psi^2\\
& -&\il_{\Si_t}(1+ \frac {t^2}{(t-u)^2})(t-u)\Theta \psi^2
\eeaa
Therefore it suffices  to show that
\be{suffices}
\il_{\Si_t}(1+ \frac {t^2}{(t-u)^2})(t-u)\Theta \psi^2\le \la^{-\eps_0}
\il_{\Si_t} 
 t^2 |\nabb \psi |^2 + (1+\frac {t^2}{(t-u)^2})\psi^2
\end{equation}
Consider the worst  term
\be{tuTheta}
\il_{\Si_t}\frac {t^2}{(t-u)}\Theta \psi^2\les
\bigg(\il_{\Si_t} {t^2}\Theta^2\psi^2\bigg)^{{\half}}
\bigg(\il_{\Si_t}
\frac{t^2}{(t-u)^2}\psi^2\bigg)^{{\half}}
 \end{equation}
According to the estimate \eqref{polar100}
of proposition \ref{Polar},  applied to   exponent $p$ such  $2p'=q$,
$$\il_{\Si_t} {t^2}\Theta^2\psi^2\les  
t^{2-\frac{4}{q}} \sup_{ u } \|\Theta \|^2_{L^{q}(S_{t,u})}
t^2\int_{\Si_t}\big( |\nabb \psi|^2 +r^{-2} |\psi|^2\big).
$$
Or, since according to \eqref{rtu}, $c^{-1} \le \frac{r}{(t-u)}\le c$,
and with the help of the estimate \eqref{trih3} for $\Theta$ with $q>2$
sufficiently close to $2$,
$$\il_{\Si_t} {t^2}\Theta^2\psi^2\les  
\la^{-5\eps_0} 
\int_{\Si_t}\big(t^2 |\nabb \psi|^2 +\frac{t^2}{(t-u)^2} |\psi|^2\big).
$$
Thus, back to 
\eqref{tuTheta}
\be{tuThetanew}
\il_{\Si_t}\frac {t^2}{(t-u)}\Theta \psi^2\les
\la^{-2\eps_0}\int_{\Si_t}\big(t^2 |\nabb \psi|^2 +\frac{t^2}{(t-u)^2} |\psi|^2\big)
\end{equation}
as desired in the proof of \eqref{suffices}. The remaining
term
on the left hand side of
\eqref{suffices} is easier to treat.

\section{Proof of  the $L^2-L^\infty $ decay estimate; theorem \ref{L2decay}}
 In this section we rely on the Boundedness Theorem \ref{BThrm} to prove the
 crucial theorem \ref{L2decay}.

Recall that  $\EE[\psi] = \EE^i[\psi] + \EE^e[\psi]$, where
\beaa
\EE^i[\psi](t)& =&  \il_{\Si_t} \Bigl (t^2 |\pr\psi|^2 +|\psi|^2\Bigr ) (1-\zeta), 
\\
\EE^e[\psi](t)& = & \il_{\Si_t} \Bigl (\ub^2 |L \psi|^2 + \ub^2|\nabb\psi|^2) + 
u^2|\Lb \psi|^2 + |\psi|^2 \Bigr )\,\zeta 
\eeaa
with a cut-off function $\zeta$ equal to $1$ in the region $u\le \frac t2$.

\medn
{\bf Estimate for $(1-\zeta) P\psi $:}\quad 
 
Observe that since the projector $P$ is an averaging operator on the scale
of size 1 and $(1-\zeta)$ is a cut-off function with the scale of size $t\ge 1$,
we can essentially write that $(1-\zeta) P \psi \approx P(\psi (1-\zeta))$.
Thus the Bernstein inequality, followed 
by the
 fact $\| (1-\zeta)\nab \psi\big)\|_{L^2(\Si_t)}\le  t^{-1} \EE^{{\half}} [\psi](t)$
 and $|\nab\zeta|\les t^{-1}$,
 implies that
\be{intdecay}
\|P (\psi(t)) (1-\zeta) \|_{L^\infty_x}\les 
\|\nab\big(\psi (1-\zeta)\big)\|_{L^2(\Si_t)}\le t^{-1} \EE^{{\half}} [\psi](t)
\end{equation}
as desired.

\medn
{\bf Estimate for $\zeta P \psi\,$:}\quad
It clearly suffices to establish the estimate for
$P\psi (t,x)$ at any point $(t,x)$ with 
$0\le u\le \frac t4$. 
According to the Sobolev inequality \eqref{asob4}, with $p=4$, of proposition \ref{Triso}
we have 
 for any positive $\de < 1 $,
$$
\sup_{S_{t,u}} |P\psi |^2\,\les\, t^{\frac {4 \de }{4+\de}}
\bigl (\il_{S_{t,u}} (|\nabb P\psi|^2 + \frac 1{t^2} |P \psi|^2) \bigr )^
{ 1-\frac {4\de }{4 +\de }}
\Bigl [{\il_{S_{t,u}} (|\nabb P\psi|^4 + \frac 1{t^4} |P\psi|^4)}\Bigr ]
^{\frac {2\de}{4+\de}}.     
$$
Using the isoperimetric inequality \eqref{asob3}  applied to $(P\psi)^2$ and $|\nabb
P\psi|^2$,
$$
\aligned
\ &\Bigl (\il_{S_{t,u}}  |P\psi|^4\Bigr )^{{\half}} \les
 \Bigl (\il_{S_{t,u}}  |\nabb P\psi|^2\Bigr )^{{\half}}
\Bigl (\il_{S_{t,u}}  |P\psi|^2\Bigr )^{{\half}} +
\frac 1t \il_{S_{t,u}}  |P\psi|^2,\\
\ &\Bigl (\il_{S_{t,u}} |\nabb P\psi|^4\Bigr )^{{\half}} \les 
\Bigl (\il_{S_{t,u}} |\nabb^2 P\psi|^2\Bigr )^{{\half}} 
\Bigl (\il_{S_{t,u}} |\nabb P\psi|^2\Bigr )^{{\half}} + 
\frac 1t \il_{S_{t,u}} |\nabb P\psi|^2.
\endaligned
$$
In addition, making use of  the trace inequality \eqref{traceineq},  
\footnote{The tensor version of the estimate requires the covariant 
$\nabb_N$ derivative. Recall that $\nabb_N$ denotes the projection on
$S_{t,u}$ of the
covariant derivative $\nab_N$. }
$$
\il_{S_{t,u}} |f|^2 \les \Bigl (\il_{\Ext} |N (f)|^2 \Bigr )^{{\half}}
 \Bigl (\il_{\Ext} |f|^2 \Bigr )^{{\half}} + 
\frac 1t \il_{\Ext} |f|^2.
$$ 
Here, $\Ext=\Si_t\cap \{0\le u \le \frac t2\}$ and 
$N$ is the vectorfield of the unit normals to $S_{t,t-\rho}$.

Thus, setting $\varep =\frac {4\de}{4+\de}$, using the fact that $t\ge 1$, and applying
the H\"older inequality, we obtain
\begin{align}
&\sup_{S_{t,u}} |P\psi |^2\,\les\,t^{\varep } 
\Bigl (\il_{\Ext} |\nabb_N \nabb P\psi |^2 + 
|\nabb P \psi |^2 + \frac 1{t^2} \big (| N (P\psi)|^2 + |P\psi|^2\big ) 
\Bigr )^{1- \varep }\cdot {\cal I}^\varep\label{extest}\\
&{\cal I}= \il_{\Si_t} |\nabb_N \nabb^2 P\psi |^2 + |\nabb^2 P\psi |^2 +
 |\nabb_N \nabb P\psi |^2 + |\nabb P\psi |^2 + \frac 1{t^4} 
(| N (P\psi)|^2 + |P\psi|^2) .\nn
\end{align}
Note that we can always replace the outside $N$ derivative with a generic derivative
$\pr$. More precisely, $|N(f)|^2 \les \sum_i |\pr_i f|^2 $.

We make the following three observations:

\medn
1)\,\,  The  derivatives in the second factor ${\cal I}$ can be ignored
in view of the presence of the projection  $P$. Thus we  can crudely bound it  by 
${\cal I}\les \il_{\Si_t} |\psi|^2\le \EE[\psi](t)$. 
  
\medn
2)\,\,The terms 
$\,\,\frac 1{t^2} \il_{\Ext} \big (| N (P\psi)|^2 + |P\psi|^2\big )$
are easily estimated by $t^{-2} \EE[\psi](t)$.

\medn
3)\,\,It remains to handle the terms 
$$
\il_{\Ext}  |\nabb_N\nabb P\psi |^2 + 
|\nabb P\psi |^2 
$$
Consider first the integral $\il_{\Ext} |\nabb P\psi |^2 $.
Let $\zeta$ be a cut-off function of the exterior region 
$\Ext$ such that $\zeta|_{\Ext} =1 $ and $|\nab \zeta|\les t^{-1}$.
We introduce the angular  vectorfields $A_i=\zeta\,( \pr_i-<\pr_i, N>N)$. Clearly,
for any scalar function $f$, $|\nabb f|^2\approx\sum_{i=1}^3|A_i f|^2$ in the
exterior region $\Ext$.
Now write,

Thus,
\beaa
\il_{\Ext} |\nabb P\psi |^2&
\approx&\sum_{i=1}^3\il_{\Ext}|A_i P\psi|^2\\
&\les&\sum_{i=1}^3\il_{\Ext}|PA_i
\psi|^2+\sum_{i=1}^3\il_{\Ext}|[P, A_i]\psi|^2\\
&\les&\sum_{i=1}^3\il_{\Si_t}|PA_i 
\psi|^2+\mbox{error}\\
&\les&\il_{\Si_t}|\nabb \psi|^2+\mbox{Error}
\eeaa
We estimate the error term $\il_{\Ext}\sum_i|[P, A_i]\psi|^2$
with the help of the following
\begin{lemma} Consider a vectorfield $X=\sum_i X^i\pr_i$ vanishing on the 
complement of the exterior
region $\Ext$ of $\Si_t$ and $P$ the standard Littlewood-Paley projection on frequencies of size
1. Then, for arbitrary scalar functions
$f$ we have the inequality\footnote{In fact the exterior region on the
right hand side of the inequality should be somewhat  enlarged( by size one ). Since
this enlargement doe not affect our arguments we prefer to ignore it. }:
$$
\|[P, X]f\|_{L^2(\Ext)}\les \sup_{i,j}\|\pr_i X^j \|_{L^\infty(\Ext)} 
\|f\|_{L^2(\Si_t)}
$$
\label{commP1}
\end{lemma}
\begin{proof}
We postpone the proof until the end of section 8, see lemma \ref{AcommP}.
\end{proof}

We apply the above lemma to the vectorfields $A_k = \zeta\,(\de_k^j - N_k N^j)\pr_j$.
Observe that the components $A_k^j$ are bounded and $|\nab \zeta |\les t^{-1}$. 
Thus
$$
\mbox{Error} \les \bigg(t^{-2} +
\|\nab N \|^2_{L^\infty(\Ext)}\bigg)\|\psi\|^2_{L^2(\Si_t)}
$$
Recall  the expression, see \eqref{Theta'},  
$\Theta = |\trch -\frac 2r| + |\chih| +|\eta| + |\pr H|$ 
and the inequality  \eqref{Dframe} $|\nab N|\les \frac 1r + \Theta$.
Observe also  that in the exterior region $\Ext$, 
$\frac 1r \le \frac 2t$. Therefore,
$$
\mbox{Error} \les \big (t^{-1} + \|\Theta \|_{L^\infty(\Ext)}\big)^2
\|\psi\|_{L^2(\Si_t)}^2
$$
We can finally conclude that
\bea
\il_{\Ext} |\nabb P\psi |^2 &\les& \il_{\Si_t}|\nabb \psi|^2 +
\big (t^{-1} + \|\Theta \|_{L^\infty(\Ext)}\big)^2
\int_{\Si_t}|\psi|^2\nn \\ &\les & 
\big(t^{-2} + \|\Theta \|^2_{L^\infty(\Ext)}\big)\EE[\psi](t)\label{Pnab}
\end{eqnarray}
We now consider $\il_{\Ext}|\nabb_N\nabb P\psi|^2$. In
view of the simple  commutation estimates \eqref{comm}
we can write:

\beaa
\il_{\Ext}|\nabb_N\nabb P\psi|^2 &\les &
 \il_{\Ext}|\nabb  (NP\psi)|^2 + \il_{\Ext}\big(r^{-1} +\Theta\big)^2|\nab P\psi|^2\\
&\approx& \sum_{i=1}^3\il_{\Ext}|A_i  (NP\psi)|^2 + \il_{\Ext}\big(r^{-1}
+\Theta\big)^2|\nab P\psi|^2
\eeaa
Observe that
 \beaa
A_i  (NP\psi)&=&A_i\,  N^j\pr_j(P\psi)=N^j A_i\pr_j(P\psi)+
[A_i, N^j]\pr_j(P\psi)\\
&=&N P(A_i\psi)+N^j[A_i, \pr_j P]\psi+[A_i, N^j]\pr_j(P\psi)
\eeaa
Therefore, using the lemma \ref{commP1}, with $P$ replaced by $\nab P$,
as well as the estimates \eqref{Dframe}
$$
\int_{\Ext}|A_i  (NP\psi)|^2\les \int_{\Ext}|A_i  \psi|^2+
\big (t^{-1} + \|\Theta \|_{L^\infty(\Ext)}\big)^2
\int_{\Si_t}|\psi|^2\nn 
$$
and finally,

\be{PnabN}
\int_{\Ext}|\nabb_N  \nabb P\psi|^2\les 
\big (t^{-2} + \|\Theta \|_{L^\infty(\Ext)}^2\big)
\EE[\psi](t) 
\end{equation}
Substituting \eqref{Pnab}-\eqref{PnabN} back into \eqref{extest} we infer
that in the exterior region 
$$
\aligned
\sup_{\stu} |P\psi|^2 &\les t^{\varep}
\big(t^{-2} + \|\Theta \|^2_{L^\infty(\Ext)}\big)^{1-\varep}
\EE^{1-\varep}[\psi](t)\cdot  {\cal I}^{\varep}\\ &\les 
t^{\varep}
\big(t^{-2} + \|\Theta \|^2_{L^\infty(\Ext)}\big)^{1-\varep}\EE[\psi](t)
\endaligned
$$
Finally, together with the interior estimates  \eqref{intdecay} this implies that
\be{Ppsi}
\|P \psi(t)\|_{L^\infty_x} \les \bigg (\frac 1{(1+t)^{1-2\varep}} + 
t^\varep\|\Theta\|_{L^\infty(\Ext)}^{1-\varep}\bigg ) \EE^{{\half}}[\psi](t).
\end{equation}
Observe that according to \eqref{trih3} of the Asymptotics Theorem 
$\Theta$ obeys the following estimate in the exterior region:
$$
\|\Theta (t)\|_{L^\infty(\Ext)}\les 
t^{-1} \la^{-\varep_0}  + \la^{\varep}\|\pr H(t)\|_{L^\infty_x}.
$$
Define
$$
d(t) = t^{\varep} \big( \la^{\varep}\|\pr H(t)\|_{L^\infty_x}
\big)^{1-\varep}
$$
Therefore,
$$
\|P \psi(t)\|_{L^\infty_x} \les \bigg (\frac 1{(1+t)^{1-2\varep}} + 
d(t)\bigg ) \EE^{{\half}}[\psi](t).
$$
To prove the desired  $L^2 - L^\infty$ decay estimate it remains to check that for some
\footnote{We can assume that $\frac 2{1-\varep}<q< 2 + 10^{-1}\ep_0$.}  
$q>2$, 
$$
t_*^{\frac 1q} \|d\|_{L^q_{[0,t_*]}}\les 1
$$
Since $t_*\le \la^{1-4\varep_0}$ it clearly suffices 
to show  that
$\|d \|_{L^q_{[0,t_*]}}\les \la^{-\frac 12}$.
In view of the estimates, see proposition \ref{propH},
$$
\|\pr H\|_{L^2_t L^\infty_x}\les \la^{-\frac 12 -4\varep_0},
\qquad
 \|\pr H\|_{L^\infty_t L^\infty_x}\les \la^{-\frac 12 + \varep_0},
$$
we infer that
\beaa
\|d\|_{L^q_{[0,t_*]}}&\les& t_*^{\varep} \la^\varep \|\pr
H\|^{1-\varep}_{L^{q(1-\varep)}_t L^\infty_x} 
\les t_*^{\varep}\la^\varep \|\pr H\|_{L^\infty_t
L^\infty_x}^{1-\frac{2}{q}-\varep}\cdot\|\pr H\|^{\frac 2q}_{L^2_t L^\infty_x}
\les\la^{-\frac 12},
\eeaa
as desired.

\section{Proof of the reduction steps}
In this section we give precise statements and proofs for the reduction
steps discussed in section 2. Recall the equation \eqref{I3}, written
in the form  \eqref{I3'}, 
\be{inA1}
\gg^{\a\b}\pr_\a\pr_\b\phi=N(\phi,\pr \phi)
\end{equation}
where $\phi=(\gg_{\mu\nu})$, $N=N_{\mu\nu}$
and $\gg^{\a\b}=\gg^{\a\b}(\phi)$. In fact $(\gg^{\a\b}) =\phi^{-1}$. We consider solutions
$\phi$ of \eqref{inA1}  such that   the  components  of both $\phi$ and $\phi^{-1}$
are uniformly  bounded.  Moreover $\gg_{\mu\nu }$ approach  the Minkowski metric $\mm_{\mu\nu}$
at infinity according to \eqref{limit}. To avoid repeating this statement in what follows
we introduce the following notation:
\begin{definition} We say that $f\in {\bold H}^s=\Hb^s(\rrrr^3)$ if $\nab f\in H^{s-1}$,
$f$ is continuous and tends to zero  as $|x|\to \infty$. Observe that $\Hb^s$,
with $s>\frac{3}{2}$, is the closure of $C_0^\infty$ in the norm $\|\nab f\|_{H^{s-1}}$.
Given a solution $\phi$ of \eqref{inA1} we say that $\phi=(\gg_{\mu\nu})\in C([0,T];
\mm+\Hb^s)$ if, for every $t\in [0,T]$, $(\gg_{\mu\nu}(t)-\mm_{\mu\nu})\in \Hb^s(\Si_t)$ and 
$\pr_t \phi\in H^{s-1}(\Si_t)$.

\end{definition}
Throughout the section
we shall use the following notation:
\begin{definition} For any function  $f$
on $\Si_t=\rrrr^3$,  $P_\la f={\cal F}^{-1}\big(\chi(\la^{-1}\xi)\hat{f}(\xi)\big)$
with $\chi$ supported in the unit  dyadic region $\half \le |\xi|\le 2$.
Also $f=\sum_{\la} P_{\la} f$. We shall  denote by  $ f_{\le\la}=P_{\le\la}f=
\sum_{\mu\le \la} f^\mu$. We shall also use the notation $f_{<\la}=P_{< \la}f=
\sum_{\mu< 2^{-M_0}\la} f^\mu$, for a sufficiently, fixed, large constant $M_0$,
 such as $100$.
\end{definition}
\begin{remark} Observe that if $f$ is continuous, approaches a constant $c$
at infinity, i.e $\sup_{|x|=r} |f(x)-c|\to 0$ 
as $r \to \infty$, and $\nab f\in H^{s-1}$,
$s>\frac{3}{2}$, then\footnote{This can be easily
proved by a density argument.} $P_\la f\in H^{s}$. 
\label{finite}
\end{remark}

 \subsection{\,\,\,Energy estimates}

We start with the following well known statement:
\begin{proposition}[Energy estimate]
Let  $ \phi\in C([0,T]; m+\Hb^s) $ be a solution of \eqref{inA1} on the time interval
$[0,T]$ for some
$s> \frac 32$ such  that $\|\phi,\phi^{-1}\|_{L^\infty_{[0,T]} L^\infty_x}\le
\La_0$.  Then $\phi$ verifies the following energy estimate.
\begin{equation}
\|\pr \phi\|_{L^\infty_{[0,T]}\dot H^{s-1}} \le C(\|\pr\phi\|_{L^1_{[0,T]} L^\infty_x},
\La_0)
\|\pr\phi(0)\|_{\dot H^{s-1}}.
\label{ens22}
\end{equation}
\label{Ens}
\end{proposition}

\begin{remark} Throughout this section we shall 
 often ignore the dependence on  $\La_0$ and the constant $M_0$ involved 
in the definition of $P_{<\la}$. 
\end{remark}
\medn
{\bf Proof:}\quad The proof of proposition \ref{Ens} can be easily reduced
to the following lemma.

\begin{lemma}
Let $\phi$ satisfy the conditions of proposition \ref{Ens}. Then
for each dyadic
$\la\in 2^{\Bbb Z}$,
$\phi^\la=P_\la\phi
$ verifies the equation
\begin{equation}
-\pr^2_t \phi^\la + (\ns^2 \gg^{0i})_{<\la}(\phi) \pr_t\pr_i \phi^\la  +
(\ns^2 \gg^{ij})_{<\la}(\phi) \pr_i \pr_j \phi^\la = R_\la,
\label{eqla}
\end{equation}
where for any $s> 1$ and $t\in [0,T]$ the right hand-side $R_\la$ has 
Fourier support in $\{\xi:\, \frac 14\la \le |\xi|\le 4\la\}$ and obeys the
estimate
\begin{equation}
\Bigl (\sum_\la \|R_\la (t)\|^2_{\dot H^{s-1}}\Bigr )^{{\half}} \le C
\|\pr\phi (t)\|_{L^\infty_x} \cdot
\|\pr \phi (t)\|_{\dot H^{s-1}} .
\label{rla}
\end{equation}
with $C$ a  constant depending only on
$ \La_0$. Moreover $\phi^\la$ also satisfies the
equation
\be{eqla222}
\gg^{\a\b}_{<\la} \pr_\a\pr_\b \phi^\la= R_\la
\end{equation}
with a different $R_\la$ which verifies the same estimate \eqref{rla} and 
the frequency property.
\label{Eqla}
\end{lemma}

\medn
{\bf Proof of lemma \ref{Eqla}}

\medn
The proof of the lemma is based on the technique of the paradifferential calculus and
is standard\footnote{The equations discussed in the literature
 are somewhat different from the one treated here because of the
non triviality of the components $\gg^{00}$ and $\gg^{0i}$ of the metric. This adds only
minor technical complications.}. For the sake of completeness
we provide an outline of the arguments. For a  more detailed treatment see for example
\cite{Ba-Ch1} or \cite{Kl-Rod}.

Recall that $P_\la$ denotes the projection on the frequencies of size $\la$,
so that $\phi^\la = P_\la \phi$.
We write the equation $\gg^{\a\b}(\phi) \pr_\a\pr_\b \phi = N$ in the form
$-\pr^2_t \phi + \ns^2 \gg^{0i}(\phi)  \pr_t \pr_i \phi + \ns^2 \gg^{ij}(\phi)  
\pr_i \pr_j \phi = \ns^2 N$.
Then 
$$
-\pr^2_t \phi^\la + P_\la (\ns^2 \gg^{0i} (\phi)  \pr_t \pr_i \phi) +
 P_\la (\ns^2 \gg^{ij}(\phi)  \pr_i\pr_j \phi ) = P_\la (\ns^2 N).
$$
 For convenience we  introduce
\be{8.1}
G\cdot \pr^2\phi = \ns^2 \gg^{0i}(\phi)  \pr_t \pr_i \phi +
\ns^2 \gg^{ij}(\phi)  \pr_i\pr_j \phi
\end{equation}
and note that at least one of the derivatives on the right hand-side
is a spatial derivative.
Then 
$$
\aligned
& P_\la (G\cdot \pr^2\phi ) = P_\la \sum_{\mu,\nu}G^\mu\cdot \pr^2\phi^\nu =
P_\la \sum_{\mu < {\half} \nu,\nu}G^\mu\cdot \pr^2\phi^\nu  + \\
&P_\la \sum_{\nu <{\half} \mu,\mu}G^\mu\cdot \pr^2\phi^\nu +
P_\la \sum_{2^{-M_0} \nu \le \mu \le 2^{M_0}\nu,\nu}G^\mu\cdot \pr^2\phi^\nu =
E_1(\la) + E_2(\la) + E_3(\la).
\endaligned
$$
It is clear that in the case when of one frequencies $\mu$ or $\nu$ 
dominate, the projection $P_\la$ on the frequencies of size $\la$
forces the dominant frequency to be of the same size. 
We say that $\mu\sim \la$ if $\frac 14\la\le \mu  \le 4\la$.

\medn
{\bf Treatment of ${\bf E_1}$}
$$ 
E_1 =  \sum_{\mu < {\half} \la}G^\mu\cdot \pr^2\phi^\la + 
\sum_{\nu\sim \la} [P_\la, G_{<{\half}\nu} \cdot] \pr^2\phi^\nu.
$$
The first term is precisely the term to keep\footnote {Observe that
$\sum_{\mu < {\half} \la}G^\mu\cdot \nab^2\phi^\la = 
(n^2 g^{0i})_{<\la} \pr_t \pr_i \phi^\la + (n^2 g^{ij})_{<\la} \pr_i \pr_j \phi - 
\sum_{\mu=2^{-M_0-1} \la}^{2^{-M_0}\la} G^\mu \cdot \pr^2\phi^\la$ and 
the second term is of the type $E_3$} on the left hand side of the equation.
To estimate the second term 
we need to make use of the standard commutator estimate, which implies
that 
$$
\|[P_\la, G_{<{\half}\nu} ] \pr^2\phi^\nu\|_{L^2_x}\le \la^{-1} 
\|\nab G_{< {\half}\nu}\|_{L^\infty_x} \|\pr^2\phi^\nu\|_{L^2_x}\le
\la^{-1} C(\La_0) \|\nab \phi\|_{L^\infty_x} \|\pr^2\phi^\nu\|_{L^2_x}.
$$
Then, since the expression $\pr^2 \phi^\nu$ contains at least one spatial 
derivative, we obtain
$$
\aligned
\|\sum_{\nu\sim \la} [P_\la, G_{{<\half}\nu} ]
\pr^2\phi^\nu\|_{\dot{H}^{s-1}}\approx &
\la^{s-1} \|\sum_{\nu\sim \la} [P_\la, G_{{\half}\nu} \cdot] \pr^2\phi^\nu\|_{L^2_x}\\
 \les
&\la^{s-1} \sum_{\nu\sim \la}  \|\pr\phi\|_{L_x^\infty} \|\pr\phi^\nu\|_{L^2_x}\\
&\les
\sum_{\nu\sim \la} \|\pr\phi\|_{L_x^\infty} 
\|\pr\phi^\nu\|_{\dot{H}^{s-1}}.
\endaligned
$$

Squaring and summing over $\la$ we obtain the bound
$$
\bigg(\sum_\la\|\sum_{\nu\sim \la} [P_\la, G_{{\le\half}\nu} ]
\pr^2\phi^\nu\|_{\dot{H}^{s-1}}^2 \bigg)^{\f12} \les\|\pr\phi\|_{L^\infty} \|\pr
\phi\|_{\dot{H}^{s-1}}.
$$
as desired.

\medn
{\bf Treatment of ${\bf E_2}$}
$$
E_2(\la) = P_\la \sum_{\mu\sim \la}G^\mu\cdot \pr^2\phi_{< {\half} \mu} 
$$
We make use of the presence of a spatial derivative in 
$\pr^2\phi_{<{\half} \mu}$ by estimating\footnote{Observe that, in view of the remark
\ref{finite} $\|G^\mu\|_{L^2_x}$ are finite.},
\beaa
\|E_2(\la)\|_{\dot{H}^{s-1}} &\le&  \la^{s-1} \sum_{\mu\sim \la} \|G^\mu\|_{L^2_x}
\| \pr^2\phi_{< {\half} \mu}\|_{L^\infty_x}\\ &\le & 
\sum_{\mu\sim \la} \la^{s-1} \mu \|G^\mu\|_{L^2_x}
\| \pr\phi_{< {\half} \mu}\|_{L^\infty_x} \les  
\sum_{\mu\sim \la}  \|\nab G^\mu\|_{\dot H^{s-1}} \| \pr\phi\|_{L^\infty_x}.
\eeaa

Thus, squaring and summing over $\la$  we obtain
$$
\bigg(\sum_\la\|E_2(\la)\|_{\dot{H}^{s-1}}^2 \bigg)^{\f12}\les 
\|\nab G\|_{\dot H^{s-1}} \| \pr\phi\|_{L^\infty_x}.
$$
Clearly, in view of our assumptions,  $G(\phi)=\phi^{-1}$ is a smooth function of $\phi$.  
 By a standard result 
on  the composition properties of  Sobolev spaces,
\begin{equation}
\|\nab G(\phi) \|_{\dot H^{s-1}}\le C(\La_0) \|\nab \phi\|_{\dot H^{s-1}}
\label{gphi}
\end{equation}
Thus,
$$
\bigg(\sum_\la\|E_2(\la)\|_{\dot{H}^{s-1}}^2 \bigg)^{\f12} \les
\|\nab \phi\|_{\dot H^{s-1}} \| \pr\phi\|_{L^\infty_x}.
$$

\medn
{\bf Treatment of ${\bf E_3}$}

\medn
$$
E_3(\la) = P_\la \sum_{2^{-M_0}\nu\le \mu\le 2^{M_0}\nu,\,\nu\ge 2^{-M_0}\la }
G^\mu\cdot \pr^2\phi^\nu .
$$
Hence,
$$
\aligned
\|E_3\|_{\dot H^{s-1}} &\le 
\la^{s-1}\sum_{2^{-M_0}\nu\le \mu\le 2^{M_0}\nu,\,\nu\ge 2^{-M_0}\la}
\|G^\mu\|_{L^\infty_x}  \|\pr^2\phi^\nu \|_{L^2_x}\\ &\le 
\sum_{2^{-M_0}\mu\le \nu\le 2^{M_0}\mu,\,\mu\ge 2^{-M_0}\la} 
\Bigl (\frac{\la}{\nu}\Bigr )^{s-1}
\|\nab G^\mu\|_{L^\infty_x}  \cdot 
\|\pr\phi^\nu\|_{\dot{H}^{s-1}}
\les  
\sum_{\nu>2^{-M_0} \la } \Bigl (\frac{\la}{\nu}\Bigr )^{s-1}
\|\nab\phi\|_{L^\infty_x} \cdot
\|\pr\phi^\nu\|_{\dot{H}^{s-1}}.
\endaligned
$$

 To check that  the multiplicative type convolution with $\nu^{(1-s)}$ maps
$l^2\to l^2$ observe that   $\sum_{\nu>\frac 14}\nu^{1-s}<\infty$,
 for $s>1$. Thus,
$$
\bigg(\sum_\la\|E_3(\la)\|_{\dot{H}^{s-1}}^2 \bigg)^{\f12} \les
\|\nab\phi\|_{L^\infty_x}\cdot  \|\pr\phi\|_{\dot{H}^{s-1}}.
$$

\medn
It remains to treat the  term $\ns^2 N(\phi,\pr\phi)$ which  depends quadratically
on $\pr\phi$.  This is standard, it can be done
in the same way as above. This ends the proof of the estimate \eqref{rla}.
It remains to prove \eqref{eqla222}. We multiply the equation \eqref{eqla},
$(\ns^2 \gg^{\a\b})_{<\la} \pr_\a\pr_\b \phi^\la = R_\la$,
by $\ns^{-2}_{<\la}$,
$$ 
- (\ns^{-2})_{<\la} \pr^2_t \phi^\la +  
(\ns^{-2})_{<\la} (\ns^2 \gg^{0i})_{<\la} \pr_t\pr_i \phi^\la +
(\ns^{-2})_{<\la} (\ns^2 \gg^{ij})_{<\la} \pr_i\pr_j \phi^\la\\
 =(\ns^{-2})_{<\la} R_\la .
$$ 
It is easy to verify that the new right 
hand-side has   the same  properties as  $R_\la$. 
Observe also that for arbitrary smooth functions $f, g$
$$
(f g)_{<\la} = f_{<\la} g _{<\la} + P_{<2\la}\big([P_{<\la},f] g\big) +
P_{<\la} \sum_{2^{-M_0} \la\le \mu\le 2^{-M_0 +1} \la} f^\mu g_{<\la} .
$$
Applying this to $f=\ns^{-2}$ and $ g= \ns^2 \gg^{\a\b}$ with $a=0,..,3$, $\b=1,..,3$, we
obtain
$$
\aligned
\gg^{\a\b}_{<\la} \pr_\a\pr_\b \phi^\la &=   
(\ns^{-2})_{<\la} R_\la
+  P_{<2\la}\bigg([P_{<\la},\ns^{-2}] \ns^2 \gg^{\a\b}\bigg) \pr_\a\pr_\b \phi^\la
\\ &+
\sum_{\a =0,..,3, \b =1,..,3} \sum_{2^{-M_0} \la\le \mu\le 2^{-M_0 +1} \la}
P_{<\la}\bigg((\ns^{-2})^{\mu} (\ns^2 \gg^{\a\b})_{<\la}\bigg) \pr_\a\pr_\b \phi^\la 
\endaligned
$$
The commutator term on the right hand-side of the expression 
above is precisely of the type $E_1(\la)$ and can be handled similarly. 
The metric  component $\ns^{-2}$ appearing  in the second term 
contains only
frequencies $\mu\ge 2^{-M_0}\la$. This allows us  to move one spatial
derivative from 
$\pr_\a\pr_\b \phi^\la$. Hence, 
the new right hand side $R_\la$  possesses the  same properties 
as the old  $R_\la$.
\qed
\begin{remark}
In the subsequent paper we shall also need the following more general 
result concerning other dyadic projections of our equation.
\begin{lemma}
Under the assumptions of lemma \ref{Eqla} we have 
$$
\gg^{\a\b}_{<\la} \pr_\a\pr_\b\phi_{<\la} = F_\la.
$$
The function $F_\la$ obeys the estimates
$$
\|F_{\la}\|_{L^1_t L^2_x}\le C \|\pr \phi\|_{L^1_t L^\infty_x} 
\|\pr \phi\|_{L^\infty_t L^2_x} ,\qquad 
\|F_{\la}\|_{L^1_t \dot H^1}\le C \|\pr \phi\|_{L^1_t L^\infty_x} 
\|\pr \phi\|_{L^\infty_t \dot H^1_x}
$$
In addition, for any dyadic $\mu\ge 1$
$$
\gg^{\a\b}_{<\la} \pr_\a\pr_\b P_{\la\mu} \phi =F_{\la,\mu},
$$
where $F_{\la,\mu}$ verifies 
$$
\aligned
&\|F_{\la,\mu}\|_{L^1_t L^2_x}\le C (\la\mu)^{-\ga}\la^{-1} \|\pr \phi\|_{L^1_t L^\infty_x} 
\|\pr \phi\|_{L^\infty_t \dot H^{1+\ga}} ,\\
&\|F_{\la,\mu}\|_{L^1_t \dot H^1}\le C \la^{-\ga} \mu^{1-\ga}\|\pr \phi\|_{L^1_t L^\infty_x} 
\|\pr \phi\|_{L^\infty_t \dot H^{1+\ga}}
\endaligned
$$
The function $\gg=\phi^{-1}$ satisfies similar equations.
\label{paper2}
\end{lemma}
The proof of lemma \ref{paper2} proceeds in the same manner as the proof
of lemma \ref{Eqla} after applying the respective projections $P_{<\la}$ and
$P_{\la\mu}$.
\end{remark}

To finish the proof of the  proposition \ref{Ens} we 
choose a large parameter $\La$ in such a way that for any 
$\la\ge \La$ the metric $(\ns^2 \gg^{ij})_{<\la}$ is uniformly elliptic.
 This is always possible
since $P_{<\la}$ is an approximation of the identity and the
original metric $(\ns^2 \gg^{ij}) $ is uniformly elliptic in $[0,T]$.

\medn
For the values of the dyadic parameter $\la \le \La$ rewrite the
equation for $\phi^\la$ in the form
$$
-\pr^2_t\phi^\la + (\ns^2 \gg^{0i})_{<\la} \pr_t \pr_i \phi^\la  +
(\ns^2 \gg^{ij})_{<\la} \pr_i \pr_j \phi^\la  = R_\la'
$$
noting that the change of the metric introduces the error term 
of the type $E_2$.

\medn
For $\la \ge \La$ we keep the form of the equation as in lemma \ref{Eqla}
$$
- \pr^2_t\phi^\la + (\ns^2 \gg^{0i})_{<\la} \pr_t \pr_i \phi^\la  +
(\ns^2 \gg^{ij})_{<\la} \pr_i \pr_j \phi^\la = R_\la
$$
In either case, the standard $H^1$ energy estimate for the wave equation
yields
$$
\|\pr \phi^\la\|_{L^\infty_{[0,T]} L^2} \le C(\La_0) (\|\pr \phi^\la(0)\|_{L^2} +  
\|R_\la\|_{L^1_{[0,T]} L^2_x}).
$$

Using lemma \ref{Eqla} and the Gronwall inequality we immediately obtain for $s>1$
$$
\|\pr \phi\|_{L^\infty_{[0,T]} \dot H^{s-1}}\les
\exp {(\|\pr\phi\|_{L^1_{[0,T]} L^\infty_x})}\|\pr \phi(0)\|_{\dot H^{s-1}}.
$$
The estimate for $s=1$ follows by standard energy
estimates without  the paradifferential
decomposition.

\subsection{Reduction to the Strichartz type estimates}
 As discussed in section 2 we need to prove the Strichartz type inequality
\eqref{Istrich100}. This is achieved by the following

\begin{theorem}[{\bf A1}]
Let $\phi\in C([0,T]; \mm+\Hb^{1+\ga})$ be a solution of \eqref{inA1} 
on the time interval $[0,T]$, $T\le 1$. 
Assume that
\begin{equation}
\|\pr\phi\|_{L^\infty_{[0,T]} H^{1+\ga}} 
+\|\pr \phi\|_{L^2_{[0,T]} L^\infty_x}\le B_0, 
\label{B0}
\end{equation}
 
 There exists a small positive exponent $\de=\de(B_0)$ such that
$\phi$ satisfies the following  local in time Strichartz type estimate,
\begin{equation}
\|\pr \phi\|_{L^2_{[0,T]}L^\infty_x}\le C(B_0)\, 
T^{\de} 
\label{st11}
\end{equation}
\label{ThrmA1}
\end{theorem}
\begin{remark} In view of the remark
 \ref{lesremark} and definition \ref{lesdefinition}
we shall treat  $B_0$ as a universal constant in what follows
 and hide the dependence on it in the notation $\les$.
\end{remark}
\subsection{The dyadic version of 
the Strichartz type estimate}
Fix a large  frequency parameter  $\La$. It easily follows from the triangle 
inequality that
for $p\in [1,\infty]$,
$$
\|\pr\phi\|_{L^p_x} \le \|\pr\phi_{le \La}\|_{L^p_x} + 
\sum_{\la> \La} \|\pr\phi^\la\|_{L^p_x}.
$$
Thus,  theorem  \ref{ThrmA1} follows from the following dyadic version of the Strichartz
type estimates for $\phi^\la=P_\la \phi$.
\begin{theorem}[{\bf A2}]
Let $\phi$ be as in
theorem \ref{ThrmA1}. 
There exists a small positive exponent $\de=\de(B_0)$ such that 
for each $\la\ge \La$, the  function 
$\phi^\la$ satisfies the Strichartz type estimate
\begin{equation}
\|\pr \phi^\la\|_{L^2_{[0,T]}L^\infty_x}\les\,c_\la 
T^{\de} 
\label{st21}
\end{equation}
with  constants $c_\la$ such that $\sum_\la c_\la\le 1$.
\label{ThrmA2}
\end{theorem}

\medn
\begin{remark} \label{bfA2} The corresponding
estimate for small frequencies, i.e. for    $\phi_{<\la}$,  follows trivially from the 
Sobolev inequality,
$$
\|\pr \phi_{<\la}\|_{L^2_{[0,T]}L^\infty_x}\les\,T^{\half} 
\|\pr \phi_{<\la}\|_{L^\infty_{[0,T]} H^{\frac 32+\ga}}\les 
\,\La^{\frac 12}  T^{\half} \|\pr\phi\|_{L^\infty_{[0,T]} H^{1+\ga}}\les
\La^{\frac 12} T^{\half}.
$$
 Since $\La$ is a fixed
large parameter, which could depend only upon $B_0$, we have the desired bound
for the low frequency part of $\phi$.
\end{remark}

\begin{remark}
We shall need the following version of the estimate 
\eqref{rla} for $R_\la$ and any $s<2+\ga$:
\begin{equation}
\|R_\la(t)\|_{\dot H^{s-1}}\les c_\la\, \|\pr\phi\|_{L^\infty_x}
\|\pr \phi\|_{H^{1+\ga}}
\label{eist233}
\end{equation}
with constants  $c_\la$: $\,\sum_\la c_\la\le 1$. The estimate \eqref{eist233}
can be easily obtained from \eqref{rla} by making use of the fact that 
the Fourier support of $R_\la$ is localized on the set 
$\{\xi:\,\,\la\le |\xi|\le 4\la\}$.
As a consequence, using the bootstrap assumption \eqref{B0}, 
we also have the estimate
\begin{equation}
\|R_\la(t)\|_{L^1_{[0,T]} \dot{H}^{s-1}}\les c_\la\,T^{\f12}
\|\pr\phi\|_{L^2_{[0,T]}L^\infty_x}
\|\pr \phi\|_{L^\infty_{[0,T]} H^{1+\ga}}\les  c_\la
\label{eist23}
\end{equation}
\end{remark}

\subsection{ Dyadic linearization and time restriction}
This step reduces  theorem \ref{ThrmA2} to a  Strichartz type estimate for the linearized
equation $\gg^{\a\b}_{<\la} \pr_\a\pr_\b \psi =0$ on smaller subintervals 
of $[0,T]$. We partition $[0,T]$ by   the intervals $I_k= [t_k, t_{k+1}],
\, \,k=0,..,\la^{8\eps_0}$ with the properties 
$|I_k|\le T \la^{-8\eps_0}$ and
$\|\pr\phi\|_{L^2_{I_k} L^\infty_x} \le \la^{-4\eps_0} B_0$. The existence of
such partition is insured by the bootstrap condition \eqref{B0}.
\begin{theorem}[{\bf A3}]
Fix $\la\ge \La$ and $k\in{\Bbb Z}\cap [0, \la^{8\eps_0 }]$ and
let $\psi$ be a solution of the linear wave equation 
$$\gg^{\a\b}_{<\la} \pr_\a\pr_\b \psi = 0$$
on the interval $I_k=[t_k, t_{k+1}]$,
verifying,  
\begin{equation}
 (2^{-10} \la)^m \|\pr \psi(t_k)\|_{L^2_x}\le \|\nab^m \pr\psi(t_k)\|_{L^2_x}\le 
(2^{10} \la)^m \|\pr \psi(t_k)\|_{L^2_x}
\label{sup1}
\end{equation}
for every $m\ge 0$.
Then there exists a sufficiently small exponent $\de>0$ such that:
\begin{equation}
\|P_\la\,\pr \psi\|_{L^2_{I_k} L^\infty_x}\les\,  
|I_k|^\de  \|\pr\psi(t_k)\|_{\dot H^{1+\de}}
\label{st31}
\end{equation}
The size of $\de$ depends only on $\eps_0, B_0$. In particular, for any 
 $\eps_0>0$, we can chose   $\de$ such  that, 
$\de<10^{-1} \ga$.
\label{ThrmA3}
\end{theorem} 

\begin{remark}
The condition \eqref{sup1} implies that, modulo a negligible
``tail'', the Fourier support of $\pr\psi(t_k)$ belongs to the set
$\{\xi:\,\, 2^{-10}\la\le |\xi|\le 2^{10}\la\}$. In general, we shall say that function
$f$ obeys the property \eqref{sup2}$_M$ if 
\begin{equation}
(2^{-M}\la)^m \|f\|_{L^2_x}\le \|\nab^m f\|_{L^2_x}\le (2^M \la)^m \|f\|_{L^2_x}
\label{sup2}
\end{equation}
\end{remark}
\medn

\begin{lemma}
\smallskip\noindent

\begin{enumerate}
\item
Assume $f$  in $\rrrr^3$ is a function whose frequency is
localized to the region $|\xi|\le 2^{-M_0}\la$ and $c\le f\le c^{-1}$
for some positive number $c$. Then $u=f^{-1}$
verifies,
\be{nabmu}\|\nab^m u\|_{L^\infty}\les (2^{-M_0}\la)^m. 
\end{equation}
\item
Assume\footnote{Recall that $M_0$ is a large positive constant} 
that $u $  verifies \eqref{nabmu} and $c\le u\le c^{-1}$. 
Let $v$ be another  function
verifying the condition \eqref{sup2}$_5$. Then\footnote{This property
is analogous to the standard paraproduct rule concerning
the multiplication of functions  $u,v$ where  the frequency
 of $v$ dominates.   } $u\cdot v$ verifies
\eqref{sup2}$_{10}$.
\end{enumerate}

\label{invprod}
\end{lemma}
\begin{proof}
The proof of 1. is based on the trivial identity 
$f\cdot f^{-1}=1$. Differentiating it 
and applying  the Leibnitz rule we conclude that, although
the Fourier support of $f^{-1}$ does 
not belong to the set
$\{\xi:\,\,|\xi|\le 2^{-M_0}\la\}$,  we still have the property,
$$
\| \nab^m (f^{-1})\|_{L^\infty_x}\les (2^{-M_0}\la)^m.
$$
The proof of 2. is once again an exercise in Leibnitz rule.
In particular, for $m=1$ we have 
$$
\aligned
\|\nab (u\cdot v)\|_{L^2_x}&\les \|\nab u\|_{L^\infty_x} \|v\|_{L^2_x} + 
  \|u\|_{L^\infty_x} \|\nab v\|_{L^2_x}\\ &\les 2^{-M_0} \la \|v\|_{L^2_x} +
2^{5}\la \|v\|_{L^2_x}\les 2^{10}\la \|u\cdot v\|_{L^2_x}
\endaligned
$$
On the other hand,
$$
\aligned
\|\nab (u\cdot v)\|_{L^2_x}&\ges  \|u\|_{L^\infty_x} \|\nab v\|_{L^2_x} -
\|\nab u\|_{L^\infty_x} \|v\|_{L^2_x}\\ &\ges 2^{-5}\la \|v\|_{L^2_x}-
2^{-M_0} \la \|v\|_{L^2_x} \ges 2^{-10}\la \|u\cdot v\|_{L^2_x}
\endaligned
$$
\end{proof}

\medn
{\bf Proof of the implication Theorem (A3) $\to$ Theorem (A2):}\quad
We shall first prove an inhomogeneous version of the Strichartz estimate \eqref{st31}
 for solutions of the equation $ \gg^{\a\b}_{<\la} \psi = F$, with the right 
hand side $F$ verifying \eqref{sup2}$_5$. Recall that
 $\gg^{\a\b}_{<\la} =P_{\le 2^{-M_0}\la}\gg^{\a\b}$. 
The Duhamel formula on the interval $I_k$ for the inhomogeneous equation 
$ \gg^{\a\b}_{<\la} \pr_\a\pr_\b\psi = F$ takes the form 
\be{soloperator}
\psi (t) = [W(t,0)] \psi[t_k] + \il_0^t W(t,s) \Bigl (
(\gg^{00}_{<\la})^{-1} F(s)\Bigr )\,ds .
\end{equation}
 with $\psi[t]$ denoting the vector $\big(\psi(t), \pr_t\psi(t)\big)$. Here $[W(t,s)]$ is the
solution operator of the homogeneous equation acting on the pair of initial data $(w_0,w_1)$
at time
$s$, and
$W(t,s)$ 
 is a solution operator corresponding to the special type of the initial
data $(0,w_1)$. We need to check that $(\gg^{00}_{<\la})^{-1} F(s)$
verifies the same conditions \eqref{sup2} as $F$.

Recall  $-\gg^{00} = \ns^{-2}$. Since $F$ verifies \eqref{sup2}$_5$, using 1. and 2. 
of lemma \ref{invprod}, we conclude that
$\big[(\ns^{-2})_{<\la}\big]^{-1} F$ verifies \eqref{sup2}$_{10}$.

We now apply theorem \ref{ThrmA3} to \eqref{soloperator}, assuming
also that the initial data $\pr\psi(t_k)$ verify the assumption
\eqref{sup2}$_{10}$,
\be{locinhstrich}
\|P_\la\,\pr \psi\|_{L^2_{I_k} L^\infty_x}\les\,  
|I_k|^\de  \Bigl (\|\pr\psi(t_k)\|_{\dot H^{1+\de}} + 
\|F\|_{L^1_{I_k} \dot H^{1+\de}}\Bigr ).
\end{equation}
Fix a sufficiently small $\eps_0$ such that $5\eps_0 + \de < \ga$.
Consider the $\la$-dyadic piece $\phi^\la$ of $\phi$, solution 
of the equation \eqref{inA1}, as in Theorem (A2). We know  that
 $\phi^\la$ verifies the
equation
$\gg^{\a\b}_{<\la} \pr_\a\pr_\b \phi^\la = R_\la$ on 
$[0,T]$ and 
the Fourier support
of $R_\la$ belongs to the set $\{\xi:\,\,\frac 14 \la \le |\xi|\le 4 \la\}$,
thus automatically satisfying property \eqref{sup2}$_5$.
We can therefore apply \eqref{locinhstrich} to  $\phi^\la$  on each $I_k$
to obtain:
$$
\aligned
\|\pr \phi^\la \|_{L^2_{[0,T]} L^\infty_x} &= \Bigl (\sum_{k=0}^{\la^{8\eps_0} -1} 
\|\pr \phi^\la \|^2_{L^2_{I_k} L^\infty_x}\Bigr )^{\half} \\ &\les\,  
 \left (\sum_{k=0}^{\la^{8\eps_0} -1}|I_k|^{2\de}  \Bigl 
(\|\pr \phi^\la (t_k)\|_{\dot H^{1+\de}} + 
\|R_\la\|_{L^1_{[0,T]} \dot H^{1+\de}}\Bigr )^2 \right )^{\half} \\ &\les
\,  
|T|^\de  \la ^{4\eps_0} 
 \Bigl (\|\pr\phi^\la \|_{L^\infty _{[0,T]} \dot H^{1+\de}} + 
\|R_\la\|_{L^1_{[0,T]} \dot H^{1+\de}}\Bigr )\\ &\les \,  
|T|^\de
\Bigl (\|\pr \phi^\la\|_{L^\infty _{[0,T]} H^{1+\ga}} + 
\|R_\la\|_{L^1_{[0,T]}\dot H^{1+4\eps_0 + \de}}\Bigr )\\
& \les \,  |T|^\de\,c_\la 
\endaligned
$$
The last two inequalities follow from the inequality $\de+5\eps_0 < \ga$
and the estimate \eqref{eist23}.

\subsection{Properties of the metric $\gg_{<\la}$}
Recall that   $\gg^{\mu\nu}_{<\la}=P_{\le 2^{-M_0}\la}(\gg^{\mu\nu})$ where
$\gg^{\mu\nu}$ is the inverse of the Lorentz  metric $\gg_{\mu\nu}=\phi$.
We shall use the notation $\gg_{<\la}$ to denote the inverse 
of $\gg^{\mu\nu}_{<\la}$. Observe that, in view of our assumption $\la\ge \La$,
$\gg_{<\la}$ defines a Lorentz metric in our spacetime region $[0,T]\times \rrrr^3$. 
It clearly depends on the solution $\phi$ of 
the quasilinear problem \eqref{inA1}.
 In the next proposition we state the properties of
the family $\gg_{<\la}$ which follow from the bootstrap condition \eqref{B0} 
on $\phi$. We denote by $ \rr_{\a\b}(\gg_{<\la})$ the components of 
Ricci curvature of the metric $\gg_{<\la}$.

\begin{proposition}
Let $\phi\in C([0,T]; \mm +  \Hb^{1+\ga})$ 
be a solution of \eqref{inA1} on  $[0,T]$, $T\le 1$. 
Assume that $\phi$ verifies the assumption \eqref{B0} of 
theorem \ref{ThrmA1}. Then the family of metrics 
$\gg_{<\la}$ 
obeys the following conditions on each  interval $I_k$ such 
that   $|I_k|\le T\la^{-8\eps_0}$,  and  
$\|\pr\phi\|_{L^2_{I_k} L^\infty_x} \le \la^{-4\eps_0} $ :
\begin{align}
&\|\pr^{1+m}\, \gg_{<\la}\|_{L^1_{I_k} L_x^\infty}\les 
\la^{-8\eps_0 +m},
\label{as1}\\
&\|\pr^{1+m}\, \gg_{<\la}\|_{L^2_{I_k} L_x^\infty}\les 
\la^{-4\eps_0 + m},
\label{as2}\\
&\|\pr^{1+m}\,\gg_{<\la}\|_{L^\infty_{I_k} L_x^\infty}\les 
\la^{\frac {1}2 -4\eps_0 + m},
\label{as3}\\
&\|\nab^{\frac 12+m}(\pr gg_{<\la})\|_{L^\infty_{[0,t_*]} L_x^2}\les \la^{\frac 12+m}\quad
\mbox{for}\quad  0\le m \le \frac12 +4\ep_0\label{as5}\\
&\|\nab^{\frac 12+m}(\pr^2 \gg_{<\la})\|_{L^\infty_{[0,t_*]} L_x^2}\les 
\la^{\f12 +m -4\ep_0}\quad
\mbox{for}\quad  -\frac12 + 4\ep_0 \le m \label{as5'}\\
&\|\nab^m \,\gg^{\a\b}_{<\la}
 \pr_\a\pr_\b \gg_{<\la}\|_{L^1_{I_k} L_x^\infty}\les
\la^{-8\eps_0+ m},
\label{as4}\\
&\|\nab^m (\nab^{\frac 12}\rr_{\a\b}(\gg_{<\la}))\|_{L^\infty_{I_k}L^2_x}\les \la^m,
\label{as800}\\
&\|\nab^m \rr_{\a\b}(\gg_{<\la})\|_{L^1_{I_k} L_x^\infty}\les
\la^{-8\eps_0 + m}.
\label{as6}
\end{align}
\label{Asg}
\end{proposition}
\begin{remark}  It suffices to prove the above estimates
for the inverse metric $\gg_{<\la}^{\mu\nu}=P_{<\la}(\gg^{\mu\nu})$.
This can be easily seen by Leibnitz rule and the non degeneracy of $\gg_{<\la}$.
On the other hand, due to the explicit  presence of $P_{\la}$, 
the estimates for $\gg_{<\la}^{\mu\nu}$ can be immediately reduced to $m=0$. 

To be precise, the argument above works only for 
the spatial derivatives $\nab$, since $P_{<\la}$ truncates the 
frequencies of $\gg^{\mu\nu}$ only with respect to the space variable $x$.
However, using the fact that $\gg_{\mu\nu}=\phi$ is a solution of 
the wave equation, one can recover the corresponding estimates for the
time derivatives. Let us illustrate this by proving the
 estimate\footnote{This is one of the few estimates with $m\ne 0$
which  we shall  actually use.} 
\eqref{as1} with $m=1$.  We assume that we have already proved \eqref{as1}-\eqref{as4}
for $m=0$. Then, clearly the derivatives $\nab^2 \gg_{<\la}$ and 
$\nab \pr_t \gg_{<\la}$ can be estimated with an additional factor of
$\la$. It remains to address the derivative 
$\pr^2_t\,\gg_{<\la}$.
Observe that 
$$
\gg^{00}_{<\la}\pr^2_t = \, \gg^{\a\b}_{<\la} 
\pr_\a\pr_\b
+ \sum_{\a =0,..,3, \b=1,..,3}\,  \gg^{\a\b}_{<\la} \pr_\a \pr_\b.
$$

The desired estimate follows from the  condition \eqref{as4} with $m=0$ and 
the fact that the second term in the previous formula contains at least one 
spatial derivative.
\end{remark}
In view of the above remark  we shall make
no distinction between $\gg_{<\la}$ and $\gg_{<\la}^{-1}$ in
what follows.

\medn
{\bf Proof of \eqref{as1}-\eqref{as6} for  $m=0$:} \quad
The proof of inequality \eqref{as2}  
follows immediately from the definition of  $I_k$,
since
$$
\|\pr\gg_{<\la}\|_{L^2_{I_k} L^\infty_x} \les 
\|\pr\phi\|_{L^2_{I_k} L^\infty_x}\les \la^{-4{\eps_0}} 
$$
Moreover, we have an  even  stronger estimate,
\begin{equation}
\|\pr \gg\|_{L^2_{I_k} L^\infty_x} \les
\|\pr\phi\|_{L^2_{I_k} L^\infty_x}\les \la^{-4 {\eps_0}} 
\label{estrg}
\end{equation}
The H\"older inequality yields \eqref{as1} from \eqref{as2}.

\medn
The estimates  \eqref{as3}, \eqref{as5}, and \eqref{as5'} follow by
a simple application of the Sobolev inequality, the composition 
properties of Sobolev spaces and the condition $\ga> 4\ep_0$.
\begin{equation}
\begin{aligned}
\|\pr(P_{<\la}\, \gg(\phi))\|
_{L^\infty_{I_k} L_x^\infty}\les
\|\pr(P_{<\la}\,\gg(\phi))\|_{L^\infty_{I_k} H^{\frac 32+\eps}}& \\ \les
\la^{{\half}-4\eps_0}\|\pr \phi\|_{L^\infty_{I_k} H^{1+\ga}} \les 
\la^{{\half}-4\eps_0}.&
\end{aligned}
\end{equation}

\medn
The most interesting part of the proposition are the estimates \eqref{as4}, 
\eqref{as6}. Recall that the original metric $\gg$ satisfied the Einstein
equation, $\rr_{\a\b} (\gg) =0$. In addition, since $(\gg^{\mu\nu})=\phi^{-1}$ 
and
$\gg^{\a\b}  \pr_\a\pr_\b \phi =N$, each component of     $\gg^{\mu\nu}$  
satisfies the equation which can be written schematically as
$\gg^{\a\b}  \pr_\a\pr_\b \,\gg^{\mu\nu} = |\pr \phi|^2$. Thus, 
\be{useless1}\|\gg^{\a\b}  \pr_\a\pr_\b \,\gg\|_{L^1_{I_k} L^\infty_x}\les
\la^{-8\eps_0}.
\end{equation}

On the other hand we recall
the   expression for  $\rr_{\a\b} (\gg) $ relative to arbitrary
coordinates,  
$$
\rr_{\a\b} (\gg) = {\half} \,\gg^{\mu\nu} (\pr^2_{\mu\b}\,\gg_{\a\nu} + 
\pr^2_{\a\nu}\,\gg_{\mu\b} - \pr^2_{\a\b}\,\gg_{\mu\nu} - 
\pr^2_{\mu\nu}\,\gg_{\a\b}) +
\,\gg_{\ga\de} (\Ga^{\ga}_{\mu\b}\Ga^{\de}_{\a\nu}-
\Ga^{\ga}_{\mu\nu}\Ga^{\de}_{\a\b}).
$$
Here $\Ga^{\ga}_{\mu\b}$ are the Christoffel symbols of the metric 
$\gg$. It is then easy to see that the equation $\rr_{\a\b} (\gg) =0$
also implies that 
\be{useless2}
\|\gg^{\mu\nu} (\pr^2_{\mu\b}\,\gg_{\a\nu} + 
\pr^2_{\a\nu}\,\gg_{\mu\b} - \pr^2_{\a\b}\,\gg_{\mu\nu} - 
\pr^2_{\mu\nu}\,\gg_{\a\b})\|_{L^1_{I_k} L^\infty_x}\les 
\|\pr \gg\|^2_{L^2_{I_k}L^\infty_x}\les \la^{-8\eps_0}.
\end{equation}
and 
\be{useless200}
\|\gg^{\mu\nu} (\pr^2_{\mu\b}\,\gg_{\a\nu} + 
\pr^2_{\a\nu}\,\gg_{\mu\b} - \pr^2_{\a\b}\,\gg_{\mu\nu} - 
\pr^2_{\mu\nu}\,\gg_{\a\b})\|_{L^\infty_{I_k}\dot H^{\frac 12}}\les 
\|\pr \gg\cdot \pr \gg\|_{L^\infty_{I_k}\dot H^{\frac 12}}\les 1.
\end{equation}
The last inequality follows from the generalized Leibnitz rule 
and the fact that $\pr\gg\in H^{1+\ga}$.

To derive the desired
estimates \eqref{as4}-\eqref{as6} 
we simply\footnote{The estimates \eqref{as800} and \eqref{as6} also 
require the following obvious estimates,
$$
\|\pr \gg_{<\la}\|^2_{L^2_{I_k}L^\infty_x}\les \la^{-8\eps_0},
\qquad 
\|\pr \gg<{\la}\cdot \pr \gg_{<\la}\|_{L^\infty_{I_k}\dot H^{\frac 12}}\les 1.
$$ }
need to apply the following
lemma to the estimates \eqref{useless1} and \eqref{useless2}. 
\begin{lemma}
 Let ${\bf A} =
 (A^{\a\b\mu\nu}_{\ga\de})$
be a fixed constant tensor. Denote $ \gg \cdot {\bf A} \cdot \pr^2 \,\gg = 
\gg ^{\ga\de} A^{\a\b\mu\nu}_{\ga\de} \pr_\a\pr_\b \,\gg_{\mu\nu}$.
Assume that the linear combination  $ \gg \cdot {\bf A} \cdot \pr^2 \,\gg $ 
of the second derivatives of the metric 
$\gg$ satisfies the estimate 
$\|\gg \cdot {\bf A} \cdot \pr^2 \,\gg \|_{L^1_{I_k} L^\infty_x}\le
c(B_0) \la^{-8\eps_0}$. Then the same estimate holds for the linear combination
associated with the metric $\gg_{<\la}$:
\begin{equation}
\|\gg_{<\la}\cdot {\bf A} \cdot \pr^2 \,\gg_{<\la} 
\|_{L^1_{I_k} L^\infty_x}\les \la^{-8\eps_0}, \qquad 
\|\gg_{<\la}\cdot {\bf A} \cdot \pr^2 \,\gg_{<\la} 
\|_{L^\infty_{I_k} \dot H^{\frac 12}}\les 1
\label{ling}
\end{equation}
\label{Boxg}
\end{lemma}
\begin{proof} 
Recall that $\gg_{<\la} = P_{<\la} \gg $. Clearly,
\begin{equation}
\|\,\gg_{<\la}-\,\gg\|_{L^2_{I_k} L^\infty_x}\les
\la^{-1} \|\nab\,\gg\|_{L^2_{I_k} L^\infty_x}\les
\la^{-1 -4\eps_0}.
\label{sq2}
\end{equation}
Then
\begin{equation}
\begin{aligned}
\|\Bigl (\gg_{<\la}-\,\gg\Bigr )\cdot {\bf A}  \cdot \pr^2 \,
\gg_{<\la} \|_{L^1_{I_k} L^\infty_x}\le 
\|\,\gg_{<\la}-\,\gg\|_{L^2_{I_k} L^\infty_x}
\|\pr^2 \,\gg_{<\la}\|
_{L^2_{I_k} L^\infty_x} & \\ \les
\la^{-1- 4\eps_0}\la  \|\pr\,
\gg_{<\la}\|_{L^2_{I_k}L^\infty_x}\les 
\la^{-8\eps_0}  &.
\end{aligned}
\end{equation}
We can now consider the term 
$\,\gg\cdot {\bf A} \cdot \pr^2 \,\gg_{<\la}$.
We have
$$
\gg \cdot  {\bf A} \cdot \pr^2 \,P_{<\la} \gg =
\gg P_{<\la} \pr  \cdot  {\bf A} \cdot \pr \,
\gg  = 
[\,\gg,P_{<\la} \pr ]\cdot  {\bf A} \cdot \pr \,
\gg  +
P_{<\la} \Bigl (\,\gg\cdot  {\bf A} \cdot \pr^2\gg+ 
\,\pr \gg\cdot  {\bf A} \cdot \pr \,
\gg\Bigr )\,
 .
$$
The commutator term can be estimated  
$$\|\Bigl ([\,\gg,P_{<\la} \pr ]\Bigr )f\|_{L^2_{I_k} L^\infty_x}\les 
\|\pr \gg\|_{L^2_{I_k} L^\infty_x}\|f\|_{L_x^\infty}\les\la^{-4\eps_0}
\|f\|_{L_x^\infty}.$$ It then follows that 
$$
\|\Bigl ([\,\gg,P_{<\la}\pr]\Bigr )\cdot  {\bf A} \cdot \pr \,
\gg  \|
_{L^1_{I_k} L^\infty_x}\les  \la^{-8\eps_0}.
$$
The remaining term satisfies the desired estimate by the assumptions of the
lemma. The proof of the $\dot H^{\frac 12}$ estimate in \eqref{ling} is similar.
\end{proof}

\subsection{Rescaling}
According to theorem \ref{ThrmA3} we need to prove a Strichartz estimate 
for any solution of the problem 
$\gg^{\a\b}_{<\la} \pr_\a\pr_\b \psi=0$ on the interval 
$I_k=[t_k,t_{k+1}]$, with initial data $\psi[t_k]=(\psi(t_k), 
\pr_t\psi(t_k))$
obeying condition
\eqref{sup1}, uniformly in $\la, k$. 

It is  convenient to replace the  above problem  by its rescaled version,
so that the initial data satisfies condition \eqref{sup1} with $\la=1$
 and the rescaled time  interval $I$ has length $\le \la^{1-8\eps_0}$.

Introduce the family of the rescaled metrics\footnote{Just 
as for $g_{<\la}$ we make no distinction between  $H_{(\la)}$,
 as Lorentz metric and its inverse.}
\begin{equation}
H_{(\la)}(t,x)=g_{<\la}(\la^{-1} (t-t_k) ,\la^{-1} x)
\label{h}
\end{equation}
We  decompose the Lorentz metric $H=H_{(\la)}$
relative to our spacetime  coordinates;
\be{decompH}
-n^2dt^2+h_{ij}(dx^i+v^i dt)\otimes(dx^j+v^j dt)
\end{equation}
where $n$ and $v$ are related to ${\bf n}$, ${\bf v}$
according to the rule \eqref{h}. In view of our choice of $\la\ge \La$
 and \eqref{timelike} it easily follows that $H$ is indeed a Lorentz
 metric and
\be{psith}
c|\xi|^2\le h_{ij}\xi^i\xi^j \le c^{-1}|\xi|^2,\quad  n^2-|v|^2_h\ge c>0, \quad
|n|,|v|\le c^{-1}
\end{equation}

Proposition \ref{Asg} implies that 
$H=H_{(\la)}$ obeys the following estimates 
on the time interval
$I=[0,t_*]$ with $t_* \le\la^{1-8\eps_0}$:

\noindent {\sl Background Estimates(see proposition \ref{propH})}:
\begin{align}
&\|\pr^{1+m}H\|_{L^1_{[0,t_*]} L_x^\infty}\les \la^{-8\eps_0},
\label{ash1}\\
&\|\pr^{1+m}H\|_{L^2_{[0,t_*]} L_x^\infty}\les
\la^{-\frac 12 -4\eps_0},
\label{ash2}\\
&\|\pr^{1+m} H\|_{L^\infty_{[0,t_*]} L_x^\infty}\les 
\la^{-\frac 12-4 \eps_0},
\label{ash3}\\
&\|\nab^{\frac 12+m}(\pr H)\|_{L^\infty_{[0,t_*]} L_x^2}\les \la^{-m}\quad
\mbox{for}\quad  -\f12\le m \le \frac12 +4\ep_0\label{ash5}\\
&\|\nab^{\frac 12+m}(\pr^2 H)\|_{L^\infty_{[0,t_*]} L_x^2}\les 
\la^{-\f12 -4\ep_0}\quad
\mbox{for}\quad  -\frac12 + 4\ep_0 \le m  \label{ash5'}\\
&\|\pr^m \big({H^{\a\b}}\pr_\a\pr_\b H\big)\|_{L^1_{[0,t_*]} L_x^\infty}\les 
\la^{-1-8\eps_0},
\label{ash4}\\
&\|\nab^m \big(\nab^{\f12}\ric(H)\big)\|_{L^\infty_{[0,t_*]} L_x^2}\les 
\la^{-1},
\label{ash800}\\
&\|\pr^m \rr_{\a\b}(H)\|_{L^1_{[0,t_*]} L_x^\infty}\les 
\la^{-1- 8\eps_0 }.
\label{ash6}
\end{align} 
 
We now formulate the rescaled version of the desired Strichartz estimate.
\begin{theorem}[{\bf A4}]
Let $\psi$ be a solution of the linear wave equation
\begin{equation}
H^{\a\b}\pr_\a\pr_\b\psi = 0,
\label{st51}
\end{equation}
on the time interval $[0,t_*]$ with $t_*\le \la^{1-8\eps_0}$.
Assume that the parameter $\la\ge \La$ for a sufficiently large
constant $\La$ and that the  metric $H$ verifies \eqref{ash1}-\eqref{ash6} 
with a sufficiently small $\eps_0 >0$.
Let $P$ be the operator of projection on the set $\{\xi:\,\,1\le |\xi|\le 2\}$ 
in Fourier space. Then there exists a small constant $\de=\de(\eps_0) >0$
such that 
\begin{equation}
\|P\, \pr \psi\|_{L^2_{[0,t_*]} L^\infty_x}\les \,  
|t_*|^\de \|\pr\psi(0)\|_{L^2_x} 
\label{st52}
\end{equation}
\label{ThrmA5}
\end{theorem} 

\medn
{\bf Remark:}\quad {\it Note that Theorem ({\bf A4}) does not
contain any assumptions on the Fourier support of the initial data $\psi[0]$.}

\subsection{Decay estimates}
A variation of the standard  $TT^*$ type argument, see \cite{Kl}, allows us to reduce the
Strichartz estimate \eqref{st52} to  a corresponding dispersive inequality, see
\eqref{st61}. In the process we replace\footnote{The two wave operators differ only
by  lower order terms in so far as the Strichartz estimates are concerned.}  the  equation
$H^{\a\b}\pr_\a\pr_\b\psi =0$  by the {\it geometric} 
wave equation  $\,\square_{H} \psi =  \frac 1{\sqrt{|H |}}
\pr_\a (H^{\a\b}\sqrt {|H |}\, \pr_\b \psi )  = 0$.
\begin{theorem}[{\bf A5}]
Let $\psi$ be a solution of the linear wave equation
\begin{equation}
\begin{split}
&{\square}_{H}\psi = 0,\\
&\psi |_{t_0}=\psi_0, \quad \pr_t \psi |_{t_0}=\psi_1
\end{split}
\label{st61}
\end{equation}
on the time interval $[0,t_*]$ with $t_*\le \la^{1-4\eps_0}$ and
with initial data $\psi[t_0]=(\psi(t_0), \pr_t\psi(t_0))$. 
We consider only large values of the parameter $\la\ge \La$. 
Assume that the metric $H$ verifies \eqref{ash1}-\eqref{ash6}.
Then there exists a function $d(t)$ obeying the condition 
\be{t*d}
t_*^{\frac 1q} \|d\|_{L^q_{[0,t_*]}}\le 1, \quad \mbox{ for some $q>2$
sufficiently close to 2},
\end{equation}
such that for all $t_0\le t\le t_*$, a fixed  arbitrary small $\eps>0$, 
and a sufficiently large integer $m$,
\begin{equation}
\|P\, \pr \psi (t)\|_{L^\infty_x}\les\,  
\bigg (\frac 1{(1+|t-t_0|)^{1-\eps}} + d(t)\bigg)
\sum_{k=0}^m \|\nab^k\psi[t_0]\|_{L^1_x}.
\label{st62}
\end{equation}
\label{ThrmA6}
\end{theorem} 
We make the final reduction by decomposing the initial data 
$\psi[t_0]$ in the physical space into a sum
of functions with essentially  disjoint supports contained in
  balls of radius ${\half}$. 
Using the additivity of the $L^1$ norm and the  standard Sobolev inequality
we can reduce the dispersive inequality \eqref{st62} to an $L^2-L^\infty$
 decay estimate.
\begin{theorem}[$L^2-L^\infty$ decay]
Let $\psi$ be a solution of the linear wave equation \eqref{st61}
on the time interval $[0,t_*]$ with $t_*\le \la^\eps_0$ and with initial data
$\psi[t_0]$ supported in the ball $B_{\half} (0)$ of radius ${\half}$
centered at the origin in the physical space.
We fix a big constant $\La$ and consider only large values of the
parameter $\la\ge \La$.
Assume that the metric $H$ verifies \eqref{ash1}-\eqref{ash6}.
Then there exists a function $d(t)$ obeying the condition \eqref{t*d}
such that for all $t_0\le t\le t_*$, an arbitrary  small $\eps>0$, 
and a sufficiently large
integer $m>0$, 
\begin{equation}
\|P\, \pr \psi (t)\|_{L^\infty_x}\les \,  
\bigg (\frac 1{(1+|t-t_0|)^{1-\eps}} + d(t)\bigg)
\sum_{k=0}^m \|\nab^k \psi[t_0]\|_{L^2_x}.
\label{B2}
\end{equation}
\label{ThrmB}
\end{theorem}

\subsection {Proof of the implication Theorem ({\bf A5}) $\to$ Theorem
({\bf A4}); Decay $\to$ Strichartz}
On this step of the reduction we assume that the family of metrics $H=H_{(\la)}$
satisfies conditions \eqref{ash1}-\eqref{ash6} and that any solution of 
the geometric wave equation $\square_{H}\psi =0$ obeys
the decay estimate $$
\|P\,\pr\psi(t)\|_{L^\infty_x}\les \bigg (\frac 1{(1+|t-t_0|)^{1-\eps}} + 
d(t)\bigg) 
\sum_{k=0}^m\|\nab^k\psi[t_0]\|_{L^1_x}.
$$ 
We need to show that under
these assumptions any solution\footnote{Remark
that we don't require any assumptions
on the initial data. This is due to the presence of the 
projection $P$ in the estimate.} of the wave equation  
${H}^{\a\b}\pr_\a\pr_\b\phi =0$ satisfies
the Strichartz estimate
$\|P\,\pr\phi\|_{L^2_{[0,t_*]} L^\infty_x}\les |t_*|^\de\|\psi[0]\|_{L^2_x}$.

First, observe that it suffices to prove the following estimate:
\begin{equation}
\|P\,\pr\phi\|_{L^q_{[0,t_*]} L^\infty_x}\les  \|\phi[0]\|_{L^2_x}
\label{HStr0}
\end{equation}
with $\de = 1-\frac 2q>0$ arbitrarily small.
Observe also that the solutions of either the geometric wave equation 
$\square_{H} \psi = F$ or the equation $H^{\a\b}_\la \pr_\a\pr_\b\psi = F$ obey
the following energy inequality for any $t, t_0\in [t_0,t_*]$:
\begin{equation}
\begin{split}
\|\pr\psi(t)\|_{L^2_x} \le  \exp (C\,\|\pr H\|_{L^1_{[0,t_*]} L^\infty_x}) 
\Bigl (\|\pr\psi(t_0)\|_{L^2_x} +
\|F\|_{L^1_{[0,t_*]} L^2_x} \Bigr ) & \\\le 2 \Bigl (\|\pr\psi(t_0)\|_{L^2_x} +
\|F\|_{L^1_{[0,t_*]} L^2_x} \Bigr ),&
\end{split}
\label{stren}
\end{equation}
where the last inequality follows \footnote{Recall that we consider 
$\la\ge \La$ for a sufficiently large constant $\La$}
 from the condition \eqref{ash1} on the metric
$H$.

Furthermore, since
$$
\square_{H} = {H}^{\a\b}\pr_\a\pr_\b 
 + \frac{1}{\sqrt{|{H}|}}\pr_\a(\sqrt{|H|} H^{\a\b}\,)\pr_\b,
$$
it is easy to show \footnote{By the Duhamel Principle we would
obtain 
$$
\|P\,\pr\phi\|_{L^q_{[0,t_*]} L^\infty_x}\le M  (\|\phi[0]\|_{L^2_x} +
\|\pr H\|_{L^1_{[0,t_*]} L^\infty_x} \|\pr\phi\|_{L^\infty_{[0,t_*]} L^2_x})
$$
and the condition \eqref{ash1} together with the energy inequality for 
$\phi$ would imply \eqref{HStr0}. }
that it suffices to establish \eqref{HStr0} for a solution of the geometric 
wave equation. We shall now prove a stronger result.
\begin{proposition}
Let $\phi$ verifies the wave equation $\square_{H} \phi =0$. Assume 
that the metric $H$ is Lorentzian \footnote{for simplicity we can assume
that the ellipticity constant of the restrictions of the metric $H$ to 
the time slices $\Si_t$ is $2$} and satisfies the condition
\begin{equation}
C\|\pr H\|_{L^1_{[0,t_*]} L^\infty_x} \le {\half} 
\label{half}
\end{equation}
for some sufficiently large positive constant $C$. We also
assume that  the conclusions
of Theorem ({\bf A5}) hold true. Then, for any $\,q>2$,
\begin{equation}
\|P\,\pr\phi\|_{L^q_{[0,t_*]} L^\infty_x}\les  \|\pr \phi(0)\|_{L^2_x},
\label{HStr}
\end{equation}

\label{Proph}
\end{proposition}

\begin{proof} As in \cite{Kl}, \cite{Kl-Rod} we start by observing that
our desired estimate  
 \begin{equation}
\|P\,\pr\phi\|_{L^q_{[0,t_*]} L^\infty_x}\le M  \|\pr \phi(0)\|_{L^2_x},
\label{HStr100}
\end{equation}
is  trivially true with a constant $M>0$ which may depend on $\la$.
 Thus we only need   to prove that the constant 
$M$ is in fact independent of $\la$.
\begin{remark} We shall 
 first prove the estimate \eqref{HStr}
for $P\,\pr_t\phi$.
\end{remark}

\begin{definition}
\label{definitionphi}
Setting  $(w_0, w_1)\in H^1(\rr^3)\times L^2(\rr^3)$,  $w=(w_0,w_1)$   we
denote by   
$\Phi(t,s;w)$  the vector 
$(\phi,\pr_t\phi)$, where  $\phi(t,s;w)$ is  the solution 
at time  $t$ of the homogeneous equation ${\square}_{H}\phi=0$ subject
to the initial data at time $s$, $\phi(s,s;w)=w_0, 
\pr_t\phi(s,s;w)=w_1$.
\end{definition}
 By a  standard uniqueness argument \footnote{which follows
from the energy estimate \eqref{stren}, which still holds under assumption \eqref{half} on
the metric $H$} 
we can easily prove the following:
\begin{equation}
\Phi\bigg(t,s;\Phi(s,t_0; w)\bigg)=\Phi(t,t_0; w) 
\label{uniq}
\end{equation}

\begin{definition}
Denote by  $ {\cal H}$  the set of vector functions
$w=(w_0,w_1)$ with $(w_0, w_1) \in  H^1(\rr^3)\times L^2(\rr^3)$.
 The scalar product
in
$\cal H$ is defined by
$$
<w,v>= \int_{\Si_0}\bigg( - H^{00} w_1\cdot v_1+
 H^{ij}\pr_i w_0\cdot \pr_j v_0\bigg)
$$
\end{definition}
\begin{remark} Observe that the above scalar product
is positive definite. Indeed $H^{00}$ is strictly negative 
and $H^{ij}$ is positive definite. To see the last assertion
let $h_{ij}$ denote the metric induced by $H$ on $\Si_t$.
In fact the metric $H$  is given by $-n^2dt^2+h_{ij}(dx^i+v^i dt)
\otimes(dx^j+v^j dt)$. Thus $H^{ij}=h^{ij}-n^{-2} v^i v^j$.  
 Observe first\footnote{ Here
 $v_i=h_{ij} v^j$.}
 that $H^{ij} v_iv_j>c |v|_h^2. $ This follows easily from
 $n^2-|v|^2_h>0$,
 see \eqref{psith}.
On the other hand, denoting by $T_v=\{ \om/ h_{ij}\om^i v^j=0\}$ the 
orthogonal complement to $v$,  we easily check that 
$H^{ij} \om_i\om_j>c|\om|^2$. This follows from the positivity
of $h$, see \eqref{psith}. Finally  $H^{ij}\om_iv_j=0$. 
\end{remark}
 Let $X=L^q_{[0,t_*]} L_x^\infty$
and  its  dual $X'=L^{q'}_{[0,t_*]} L^1_x$. Let $\TT$ be the
operator
 from ${\cal H}$ to $X$   defined by:
\begin{equation}
\TT(w)=-P\,\pr_t\phi(t,0;w)
\label{ti}
\end{equation}
with $\phi$ defined according to definition \ref{definitionphi}.

The adjoint $\TT^*$ is defined from $X'$ to $\cal H$.
To prove  the estimate \eqref{HStr} it suffices
to check that $\TT\cdot \TT^*$ is a bounded operator from
$X'$ to $X$. In view of \eqref{HStr100}  
we have\footnote{We may assume that $M$ is the smallest constant
 for which \eqref{HStr100} holds true.} $\|\TT\|_{\cal
H\rightarrow  X}=M$ where  
$\|\TT\|_{\cal H\rightarrow X}$ denotes the operator norm of $\TT$.
  Thus, 
$$\|\TT\cdot \TT^*\|_{X'\rightarrow  X}=M^2.$$

  To calculate $\TT^*$ we write,
$$
<\TT^*f, w>:=<f, \TT(w)>=-\il_{[0,t_*]\times \rr^3} \pr_t\phi Pf dtdx\\
=\il_{[0,t_*]\times \rr^3} \pr_t\phi\,\bar{\square}_{H}\psi,
$$
where $\psi$ is the unique solution to the equation
\begin{equation}
\begin{split}
&\bar{\square}_{H}\psi=  \pr_\a (H^{\a\b} \pr_\b \psi ) = -Pf,  \\
&\phi(t_*)=\pr_t\phi(t_*)=0.
\end{split}
\label{eqpsi}
\end{equation}

Consequently, integrating by parts,  
we obtain
$$
<\TT^*f, w>= -\il_{\Si_0}\bigg(\pr_t\phi H^{0\b}\pr_\b\psi -
 H^{0\b}\pr_\b \pr_t\phi\,\psi\bigg) 
+\il_{[0,t_*]\times \rr^3}\bar{\square}_{H} \pr_t\phi
\, \psi . 
$$
Observe that 
$$
\bar{\square}_H\pr_t \phi = \pr_t\bar{\square}_H \phi 
- \pr_\a (\pr_t (H^{\a\b}) \pr_\b\phi).
$$
Therefore, integrating by parts once more, we have
$$
\aligned
<\TT^*f, w>=& -\il_{\Si_0}\bigg(\pr_t\phi H^{0\b}\pr_\b\psi -
 H^{0\b}\pr_\b \pr_t\phi\,\psi \\ &\qquad + \bar{\square}_H \phi\, \psi 
- 
\pr_t (H^{0\b}) \pr_\b\phi\, \psi \bigg) \\
&-\il_{[0,t_*]\times \rr^3}\bigg(\bar{\square}_{H} \phi
\, \pr_t \psi - \pr_t (H^{\a\b}) \pr_\b\phi\,\pr_\a \psi \bigg). 
\endaligned
$$
Further note that 
$-H^{0\b}\pr_\b \pr_t\phi+ \bar{\square}_H \phi - 
\pr_t (H^{0\b}) \pr_\b\phi = \pr_i ( H^{i\b} \pr_\b\phi)$,
and therefore,
$$
\aligned
\il_{\Si_0}\bigg(\pr_t\phi H^{0\b}\pr_\b\psi -
 H^{0\b}\pr_\b \pr_t\phi\,\psi + \bar{\square}_H \phi\, \psi  - 
\pr_t (H^{0\b}) \pr_\b\phi\, \psi \bigg)&\\ = 
\il_{\Si_0}\bigg(\pr_t\phi H^{0\b}\pr_\b\psi - 
H^{i\b} \pr_\b\phi \pr_i \psi\bigg )& \\= 
\il_{\Si_0}\bigg(H^{00} \pr_t\phi \pr_t\psi - 
H^{ij} \pr_i\phi \pr_j \psi\bigg ) &.
\endaligned
$$
Thus, since  
$$\square_H=\bar{\square}_H + H^{\a\b} \frac{\pr_\a \sqrt{|H|}}{\sqrt{|H|}}
\pr_\b$$ and $\square_H\phi=0$

\beaa
<\TT^*f, w>&=&\il_{\Si_0}\bigg(-H^{00} \pr_t\phi \pr_t\psi + 
H^{ij} \pr_i\phi \pr_j \psi\bigg )\\
&-&\il_{[0,t_*]\times \rr^3}
\bigg(\bar{\square}_{H} \phi
\, \pr_t \psi - \pr_t (H^{\a\b}) \pr_\b\phi\,\pr_\a \psi \bigg) \\
&=&\il_{\Si_0}\bigg(-H^{00} \pr_t\phi \pr_t\psi + 
H^{ij} \pr_i\phi \pr_j \psi\bigg )\\
&+&\il_{[0,t_*]\times \rr^3}
\bigg(H^{\a\b} \frac{\pr_\a \sqrt{|H|}}{\sqrt{|H|}}
\pr_\b\phi
\, \pr_t \psi + \pr_t (H^{\a\b}) \pr_\b\phi\,\pr_\a \psi \bigg) 
\eeaa
Thus, since $\phi[0]=w$ and recalling the definition of 
$<\,,\,>_{{\cal H}}$
$$
<{\cal T}^*f, w>=
<\psi[0]\,,w>+<R(f), w>
$$
 with $R(f)$ the linear  operator  from $X'$ to $\HH$ defined  by
the formula, 
\be{Roff}< R(f)\,,\,  w>= \il_{[0,t_*]\times \rr^3}
\bigg(H^{\a\b} \frac{\pr_\a \sqrt{|H|}}{\sqrt{|H|}}
\pr_\b\phi
\, \pr_t \psi + \pr_t (H^{\a\b}) \pr_\b\phi\,\pr_\a \psi \bigg) 
\end{equation}
Therefore,
\be{Tstarf}
\TT^*f = \psi[0] + R(f)
\end{equation}
with $\psi[0]=(\psi(0),\pr_t\psi(0))$. 
 
Henceforth,
\begin{equation}
\TT\TT^*f=\TT\psi[0] + \TT R(f)
\label{TT*form}
\end{equation}

Observe that 
$\square_{H}\psi=- Pf+e$ with
$e=H^{\a\b} \frac{\pr_\a \sqrt{|H|}}{\sqrt{|H|}}
\pr_\b\psi $.
Thus we can write $\psi= - \psi_1+\psi_2$ with,
$$
\aligned
&\square_{H}\psi_1=Pf,\\
&\square_{H}\psi_2=e
\endaligned
$$
with both $\psi_1,\psi_2$ verifying zero
initial conditions at $t=t_*$ as in \eqref{eqpsi}.  
Now   $\TT\psi[0]=-\TT\psi_1[0]+\TT\psi_2[0]$ 
 and, recalling the definition \ref{definitionphi}, 
$\TT\psi_1[0]= -P\pr_t\phi \big(t,0; \psi_1[0]\big)$.
According to the Duhamel principle, as in \eqref{soloperator} we have, with 
$\psi[t]=\big(\psi(t), \pr_t\psi(t)\big)$,
$$\psi_1[t]=\il_{t_*}^t\Phi(t,s; F(s))ds$$
 with $F(s)=(0,(H^{00})^{-1}Pf(s))=
(0, -n^{-2}Pf(s))$   and 
therefore,
$$\psi_1[0]=-\il_0^{t_*}\Phi(0,s; F(s))ds$$ 
and, in view of \eqref{uniq},
$$
\TT\psi_1[0]= P\pr_t\phi\bigg(t,0; \il_0^{t_*}
\Phi(0,s;F(s))ds\bigg)= P\il_0^{t_*}  \pr_t\phi(t,s;F(s))ds.
$$
We are now in a position to apply the dispersive 
inequality of Theorem ({\bf A6}).  
$$\|P\,\pr_t\phi(t,s;F(s))\|_{L^\infty}\le
  C\bigg((1+|t-s|)^{-1+\eps}+d(t)\bigg)\sum_{k=0}^m \|\nab^k (n^{-2}Pf(s))\|_{L^1}. 
$$
In view of \eqref{psith} and \eqref{ash3}, we have  $\|\nab^k n^{-2}
\|_{L^\infty}\les 1$. Thus, since
$P$ is the  projection on the frequencies of size $1$, we infer that
$$
\|P\,\pr_t\phi(t,s;F(s))\|_{L^\infty}\le
C\bigg((1+|t-s|)^{-1+\eps}+d(t)\bigg)\|f(s))\|_{L^1}. 
$$
Therefore, by the Hardy-Littlewood-Sobolev inequality,
\begin{equation}
\|\TT\psi_1[0]\|_{L^q_{[0,t_*]} L_x^\infty}\le C
\|f\|_{L_{[0,t_*]}^{q'}L_x^1} + 
\|\int_0^{t_*}d(t)\|f(s)\|_{L^1}\,ds\|_{L^q_{[0,t_*]}}\nn
\end{equation}
We can now make use of the assumption  \eqref{t*d} of
Theorem (${\bf A5}$) and infer that,
$$
\|\int_0^{t_*}d(t)\|Pf(s)\|_{L^1}\,ds\|_{L^q_{[0,t_*]}}\le
C t_*^{\frac 1q} \|d\|_{L^q_{[0,t_*]}}\|f\|_{L_{[0,t_*]}^{q'}L_x^1} 
\le C \|f\|_{L_{[0,t_*]}^{q'}L_x^1} 
$$
Thus
\begin{equation}
\|\TT\psi_1[0]\|_{L^q_{[0,t_*]} L_x^\infty}\le C
\|f\|_{L_{[0,t_*]}^{q'}L_x^1}
\label{StrPsi}
\end{equation}
with $C$ a constant, independent of $\la$.

To estimate $\TT\psi_2[0] $ we apply the Strichartz
inequality with a bound $M$,  see\eqref{HStr100}, 
$$\|\TT\psi_2[0]\|_{L^q_{[0,t_*]} L_x^\infty}\le M\|\psi_2[0]\|_{\cal H}
$$
where,
$$
\|\psi_2[0]\|_{\cal H}=\sup_{  \|w\|_{\cal H}\le 1}
<w,\psi_2[0]>_{\cal H} \le C \|\pr\psi_2(0)\|_{L^2}.
$$
We shall now make use of the energy estimate \eqref{stren} 
for  $\psi_2$ verifying
the equation $\square_{H}\psi_2 = e $, subject to the initial
conditions $\psi_2(t_*)=\pr_t\psi_2(t_*)=0$,
$$
\|\pr\psi_2(0)\|_{L^2_x} \le C \|e\|_{L^1_{[0,t_*]} L^2_x}\le
C \|\pr H\|_{L^1_{[0,t_*]} L^\infty_x} \|\pr\psi\|_{L^\infty_{[0,t_*]} L^2_x}
$$
Therefore, with the help of the condition \eqref{half}, we have
\begin{equation}
\|\TT\psi_2[0]\|_{L^q_{[0,t_*]}L^\infty_x}\le 
\frac 14 M  \|\pr\psi\|_{L^\infty_{[0,t_*]} L^2_x}
\label{StrPsi2}
\end{equation}

We shall now estimate the other error term  $\TT Rf$.
Since the   operator norm of $\TT$
is bounded by $M$,
 $$ \|\TT R(f)\|_{L^q_{[0,t_*]} L^\infty_x}\le M\|R(f)\|_{\cal H}.$$
On the other hand,
$$
\aligned
\|R(f)\|_{\cal H}=&\sup_{  \|w\|_{\cal H}\le 1}
<w,R(f)>_{\cal H}\\
= &-\sup_{ \|w\|_{\cal H}\le 1}\il_{[0,t_*]\times \rr^3}
\bigg(H^{\a\b} \frac{\pr_\a \sqrt{|H|}}{\sqrt{|H|}}
\pr_\b\phi
\, \pr_t \psi + \pr_t (H^{\a\b}) \pr_\b\phi\,\pr_\a \psi \bigg) 
\endaligned
$$
Estimating
in a straightforward manner we derive,
$$
\|R(f)\|_{\cal H}\le
C\|\pr H\|_{ L^1_{[0,t_*]} L^\infty_x}
\|\pr\phi\|_{L^\infty_{[0,t_*]} L^2_x} \|\pr\psi\|_{L^\infty_{[0,t_*]} L^2_x}.
$$
We use the energy inequality \eqref{stren}
to estimate $\|\pr\phi\|_{L^\infty_{[0,t_*]}  L^2_x}$.  Since the initial data 
   $\|w\|_{\cal H}\le 1$
we infer that, 
$\|\pr\phi\|_{L^\infty_{[0,t_*]} L^2_x}\le C$.
Therefore, with the help of \eqref{half}, we have
\begin{equation}
\|\TT R(f)\|_{L^q_{[0,t_*]} L^\infty_x}\le \frac 14 M 
\|\pr \psi\|_{L^\infty_{[0,t_*]} L^2_x}.
\label{StrPsi3}
\end{equation}

To estimate $\|\pr \psi\|_{L^\infty_{[0,t_*]} L^2_x}$ 
we rely on the following:

\begin{lemma} 
The solution $\psi$
of the equation $\bar{\square}_{H}\psi=-Pf$, $\psi(t_*)=
\pr_t\psi(t_*)=0$  verifies the estimate,
\begin{equation}
\|\pr \psi\|_{L^\infty_{[0,t_*]} L^2_x}\le 2 M
\|f\|_{L^{q'}_{[0,t_*]} L^1_x }
\label{offStr}
\end{equation}
\label{LoffStr}
\end{lemma}

Gathering together \eqref{StrPsi},\eqref{StrPsi2},\eqref{StrPsi3} and \eqref{offStr}
we infer that,

$$
\| \TT \TT^* f  \|_{X} =  
 \|\TT(\psi_1[0] +\psi_2[0]  + R(f))\|_{L^q_{[0,t_*]} L^\infty_x}
\le (C+  {\half} M^2) 
\|f\|_{L^{q'}_{[0,t_*]} L^1_x}
$$
Therefore, in view of  \eqref{TT*form},
$$ M^2=\|\TT \TT^*\|_{X'\rightarrow X} \le 
(C + {\half} M^2).$$

Thus we infer that  $M$ is a universal constant,
as desired.

\vskip 2pc

It only remains to prove the lemma 
\ref{LoffStr}. We proceed as follows.
 Let $t$ be fixed
in the interval $[0,t_*]$. We rewrite  the equation
 $\square_{H}\phi=0$ in the form,
\begin{equation}
\bar {\square}_{H} \phi= F = - H^{\a\b}\frac{\pr_\a \sqrt{|H|}}{ \sqrt{|H|}}
\pr_\b\phi
\label{eqq1}
\end{equation}
with  initial data  $\phi(t)=w_0, \pr_t\phi(t)=w_1$,
and $(w_0,w_1)=w\in {\cal H}_{t}$,  $\|w\|_{{\cal H}_{t}}\le 1$.
Here, the space ${\cal H}_t$ is defined by the scalar
product $\,<w,v>_{{\cal  H}_t} = \il_{\Si_t} - H^{00}w_1\,v_1 + 
H^{ij} \pr_i w_0 \,\pr_j v_0$.
We also recall that, see \eqref{eqpsi},
\begin{equation}
\bar{\square}_{H} \psi = - P f
\label{eqq2}
\end{equation}
with initial data $\psi_1(t_*)=\pr_t\psi_1(t_*)=0$.
As in \cite{Kl} and \cite{Kl-Rod} we multiply \eqref{eqq1}
by $\pr_t \psi$ and \eqref{eqq2} by $\pr_t\phi$ after which
we sum and integrate on our spacetime slab $[t,t_*]\times \rrrr^3$.
Observe that,
$$
\pr_\a (H^{\a\b}\pr_\b \psi) = (\pr_\a H^{\a\b})\pr_\b \psi + 
H^{\a\b} \pr_\a\pr_\b \psi
$$
$$
\aligned
H^{\a\b} \pr_\a\pr_\b \psi\,\pr_t\phi  + H^{\a\b} \pr_\a\pr_\b \phi\,\pr_t\psi &=
H^{\a \b} \pr_\a (\pr_t \phi \pr_\b \psi ) 
+ H^{\a \b} \pr_\b   (\pr_t \psi \pr_\a \phi )\\ & - 
H^{\a \b} ( \pr_\a \pr_t \phi \pr_\b \psi ) - H^{\a \b} ( \pr_\b  \pr_t \psi \pr_\a \phi )\\
& = H^{\a \b} \pr_\a (\pr_t \phi \pr_\b \psi ) 
+ H^{\a \b} \pr_\b   (\pr_t \psi \pr_\a \phi )\\ & - H^{\a \b}\pr_t 
( \pr_\a \phi \pr_\b \psi ) 
\endaligned
$$
Thus
$$
\aligned
\pr_\a (H^{\a\b}\pr_\b \psi)\pr_t \phi + \pr_\a (H^{\a\b}\pr_\b \phi)\pr_t \psi &=
\pr_\a( H^{\a \b} \pr_t \phi \pr_\b \psi + \pr_t \psi \pr_\b \phi) \\ & - 
\pr_t (H^{\a \b} \pr_\a \phi \pr_\b \psi ) + (\pr_t H^{\a \b}) \pr_\a \phi \pr_\b \psi\\
& =  \pr_i( H^{i \b} \pr_t \phi \pr_\b \psi + \pr_t \psi \pr_\b \phi) \\ &
+ \pr_t ( - H^{00} \pr_t\phi \pr_t\psi + H^{ij} \pr_i \phi \pr_j \psi )\\ &
+(\pr_t H^{\a \b}) \pr_\a \phi \pr_\b \psi
\endaligned 
$$
Integrating in the region $[t,t_*]\times \rr^n$ 
 we derive the identity,
$$
 \int_{\Si_t}\bigg(- H^{00} \pr_t\phi\pr_t\psi +
H^{ij}\pr_i\phi\pr_j\psi\bigg) =
-\il_t^{t_*}\il_{\Si_\tau }\bigg(-\pr_t\phi \,Pf + \pr_t\psi\, F + 
\pr_t(H^{\a\b}\,)\pr_\a\phi\,
\pr_\b\psi\bigg) .
$$
Therefore,
$$
\|\pr\psi(t)\|_{L^2}\le \|P\,\pr_t\phi\|_{L^q_{[0,t_*]}L^\infty_x}
\|f\|_{L^{q'}_{[0,t_*]}L^1_x}
+ C\|\pr H\|_{L^1_{[0,t_*]} L^\infty_x}
\|\pr\phi\|_{L^\infty_{[0,t_*]}  L^2_x}
\|\pr\psi\|_{L^\infty_{[0,t_*]} L^2_x}
$$
We recall that according to
our assumption  $\|P\,\pr_t\phi\|_{L^q_{[0,t_*]} L^\infty_x}\le
M \|w\|_{{\cal H}_{t}}\le M$.
Also  according to the energy estimate,
$\|\pr\phi\|_{L^\infty_{[0,t_*]} L^2_x}\le 2\|w\|_{\HH_{t}}\le 2$.
Therefore,

$$
\|\pr\psi\|_{L^\infty_{[0,t_*]} L^2_x}\le  
M\|f\|_{L^{q'}_{[0,t_*]}L^1_x}
+C\|\pr H\|_{L^1_{[0,t_*]}L^\infty_x}
\|\pr\psi\|_{L^\infty_{[0,t_*]} L^2_x}
$$
and therefore, since  $ C\|\pr H\|_{L^1_{[0,t_*]}L^\infty_x}\le
\frac{1}{2}$, we conclude that,
$$
\|\pr\psi\|_{L^\infty_{[0,t_*]} L^2_x}\le
2M\|f\|_{L^{q'}_{[0,t_*]}L^1_x}
$$
as desired.

\bigskip

To prove the Strichartz estimate  for
the spatial derivatives we   rely
on the proof, given above,  for $P\,\pr_t\phi$.
We thus assume that the estimate \eqref{LoffStr}
holds true for $P\,\pr_t\phi$ with a universal constant $M$.

To estimate $\|P\,\pr_k \phi\|_{L^q_{[0,t_*]} L^\infty_x}$ it suffices
to estimate the integral,
$ {\cal I}=\il_{[0,t_*]\times \rr^3 } P\,\pr_k \phi\, f dt dx $ for
 functions $f$ with 
$\|f\|_{L^{q'}_{[0,t_*]} L^\infty_x}\le 1.$  Let $\psi$
verify the equation $\bar{\square}_{H}\psi=Pf$ with 
$\psi(t_*)=\pr_t\psi(t_*)=0$.
Integrating by parts as before we infer that
$$
\aligned
 {\cal I}=\il_{[0,t_*]\times \rr^3} \pr_k \phi\,  \bar{\square}_{H}\psi 
&=\il_{\Si_0}H^{0\b}\bigg(\pr_k \phi\,\pr_\b\psi +\pr_k \psi\,\pr_\b\phi\bigg)
\\ &-
\il_{[0,t_*]\times \rr^3}\bigg(\bar{\square}_{H}\phi \,\pr_k \psi -
(\pr_k H^{\a\b})\pr_\a\phi \pr_\b \psi\bigg)
\endaligned
$$
Once again
$$
|\il_{[0,t_*]\times \rr^3}\bar{\square}_{H}\phi\, \pr_k \psi |\le
C\|\pr H\|_{L^1_{[0,t_*]}L^\infty_x}
\|\pr\phi\|_{L^\infty_{[0,t_*]} L^2_x}\|\pr\psi\|_{L^\infty_{[0,t_*]} L^2_x}
$$
Also,
$$
\il_{\Si_0}H^{0\b} \bigg(\pr_k \phi\,\pr_\b\psi +\pr_k \psi\,\pr_\b\phi   \bigg)
\le \|\pr \phi(0)\|_{L^2}
\|\pr\psi\|_{L^\infty_{[0,t_*]} L^2_x}.
$$
The energy estimate \eqref{stren} gives
$\|\pr\phi\|_{L^\infty_{[0,t_*]} L^2_x}\le 2\|\pr \phi(0)\|_{L^2}$.
According to the lemma \ref{LoffStr} we have, 
$$\|\pr\psi\|_{L^\infty_{[0,t_*]}
L^2_x}\le 2M\|f\|_{L^{q'}_{[0,t_*]} L^1_x}.
$$
Observe that the $M$ in  lemma \ref{LoffStr} depends only
on the Strichartz estimate \eqref{HStr} for $ P\,\pr_t\phi$
which we have already proved. 
Therefore,
$$
|{\cal I}|\le  CM  \|\pr \phi(0)\|_{L^2} (1+\|\pr H\|_{L^1_{[0,t_*]}L^\infty_x})
\|f\|_{L^{q'}_{[0,t_*]} L^1_x}\le CM\|\pr \phi(0)\|_{L^2}
$$
which implies,
$ \|P\,\pr_a \phi\|_{L^q_{[0,t_*]} L^\infty_x}\le CM\|\pr \phi(0)\|_{L^2}$
as desired.
\end{proof}

\subsection{Commutator lemma}
We conclude this section by presenting the proof of lemma \ref{commP1}
from section \ref{L2decay}.
Recall that the definition of the exterior region $\Ext=\{u\le t/2\}$.
\begin{lemma} Consider a vectorfield $X=\sum_i X^i\pr_i$ vanishing on the complement 
of the exterior
region $\Ext$ of $\Si_t$ and $P$ the standard Littlewood-Paley projection on frequencies of size
1. Then, for arbitrary scalar functions
$f$ we have the inequality:
\begin{equation}
\|[P, X]f\|_{L^2(\Ext)}\les \sup_{i,j}\|\pr_i X^j \|_{L^\infty(\Ext)} 
\|f\|_{L^2(\Si_t)}
\label{XP}
\end{equation}
\label{AcommP}
\end{lemma}
\begin{proof}
First observe, by expanding $X=X^j\pr_j$ relative to our system of our
coordinates on $\Si_t$, that $[P,X]=[P\pr_j, X^j] -  P(\pr_j X^j)$.
We shall denote $P_j =   P\pr_j$, the modified cut-off of the unit frequencies.
In what follows, the roles of $P$ and $P_j$ are identical.
The convolution kernels of $P, P_j$ are represented by the smooth functions
$P(x), P_j(x)$ verifying the condition that 
$|P(x)|, |P_j(x)|\les |x|^{-k}$ for any $k>0$ and $|x|\ge 1$.
In particular, for any functions $w, v$
$$
\aligned
[P, w] v  &= \il_{\Si_t} P(x-y)\big(w(y)-w(x)\big)v(y)\,dy \\ &= 
-  \il_0^1 \il_{\Si_t} P(x-y)(x-y)^i \pr_i w(\tau x + (1-\tau)y)  
v(y)\,dy \,d\tau 
\endaligned
$$
As a consequence,
\begin{equation}
\|[P,w] v\|_{L^2(\Si_t)}\les \|\nab w\|_{L^\infty(\Ext)}
\|v\|_{L^2(\Si_t)}
\label{Pwv}
\end{equation}
Similar inequality also holds for $P_j$.

We shall show that
$$
\|[P_j, X^j] f\|_{L^2(\Si_t)}+ \|P\big((\pr_j X^j) f\big)\|_{L^2(\Si_t)} 
\les \sup_{i,j}\|\pr_i X^j \|_{L^\infty(\Ext)} 
\|f\|_{L^2(\Si_t)} 
$$
Since all $X^j$ vanish outside of $\Ext$ and $P$ is a bounded operator 
on $L^2(\Si_t)$, we can easily estimate the second term, 
$$
 \|P\big((\pr_j X^j) f\big)\|_{L^2(\Si_t)}\les 
 \sup_{i,j}\|\pr_i X^j \|_{L^\infty(\Ext)} 
\|f\|_{L^2(\Si_t)}
$$ 
According to \eqref{Pwv} we also have
$$
\|[P_j, X^j] f\|_{L^2(\Si_t)}\les \sup_{i,j}\|\pr_i X^j \|_{L^\infty(\Ext)} \|
f\|_{L^2(\Si_t)}
$$

\end{proof}

\end{document}